\pdfoutput=1
\documentclass[12pt,titlepage]{article}
\usepackage{amsthm,amssymb,amsmath}
\usepackage{latexsym}
\usepackage{float,subfigure}
\usepackage{graphicx}
\newcommand{\Z}{\mathbb{Z}}
\newcommand{\R}{\mathbb{R}}
\newcommand{\C}{\mathbb{C}}
\newcommand{\re}{\operatorname{Re}}
\newcommand{\im}{\operatorname{Im}}
\newcommand{\ext}[0]{\ensuremath{\mathrm{ext}}}
\numberwithin{equation}{section}
\numberwithin{figure}{section}

\theoremstyle{plain} %text of this environment is typesetted in italics
\newtheorem{theorem}{Theorem}[section]

\newtheorem{corollary}[theorem]{Corollary}
\newtheorem{proposition}[theorem]{Proposition}
\theoremstyle{definition} %text of this environment is typesetted in roman letters
\newtheorem{definition}[theorem]{Definition}
\newtheorem{remark}[theorem]{Remark}

\begin{document}
\begin{title}
{Triply Periodic Minimal Surfaces Bounded by Vertical Symmetry Planes}
\end{title}

\author
{Shoichi Fujimori\thanks{Partially supported by JSPS Grant-in-Aid for 
Young Scientists (Start-up) 19840035.} 
\and 
Matthias Weber\thanks{This material is based upon
work for the NSF under Award No.DMS - 0139476.}
}
%\end{author}

\date{\today}
\maketitle

\noindent {\sc Abstract. } {\footnotesize}
We give a uniform and elementary treatment of many classical and new triply periodic
minimal surfaces in Euclidean space, based on a Schwarz-Christoffel formula for periodic polygons in the plane. Our surfaces share the property that vertical symmetry planes cut them into simply connected pieces.  

\noindent
{\footnotesize %2000 MSC numbers
}

\noindent
{\footnotesize %key words and phrases
}

\parindent=0cm

\section{Introduction}

By a {\em triply periodic minimal surface} $\Sigma$ we mean a complete, embedded minimal surface in Euclidean {three-}space {$\R^3$} which 
is invariant under three linearly independent translations. The largest group $\Lambda$ of Euclidean translations that act as orientation preserving
isometries on the surface is called the {\em period lattice} of $\Sigma$. The quotient $\Sigma/\Lambda$ is a compact Riemann surface of genus $g\ge 3$.

The theory of triply periodic minimal surfaces begins around 1890 with the
work of H.~A.~Schwarz and his students (\cite{schw1,ne1}), who found several now famous examples by explicitly solving suitable Plateau problems.

The NASA scientist A.~Schoen found many more examples, using experiments, numerical evidence, and the Weierstrass representation (\cite{sch1}). His
findings were ignored by the mathematical community until
H. Karcher verified his claims (\cite{ka5}) and added even more examples to the list (\cite{kapo1}). 

Since then, a vast number of more examples have been found by many people. Two important results deal with 5-dimensional families: Meeks found an explicit such family for genus $3$ (\cite{me6}), and Traizet constructed 5-dimensional families for any genus close to degenerate limits (\cite{tr3a, tr7}). All other constructions results target 1- or 2-dimensional families, using
strong symmetry assumptions. Examples with existence proofs can be found in the work of Fischer and Koch (\cite{fk1,fk2}),  Ramos Batista (\cite{Ra1}), or Huff (\cite{huff1}).

We will  here give a new construction method for triply periodic minimal surfaces which are cut by symmetry planes
into minimal surfaces with boundary that are simply connected and invariant under a vertical translation. The Weierstrass representation of these
surfaces can be given in terms of a Schwarz-Christoffel formula for {\em periodic} Euclidean polygons, using $\vartheta$-functions on suitable tori.

In the simplest case, there is no period condition, and we obtain  surfaces known to Schwarz and Schoen. The next complicated cases impose 1-dimensional period problems, which we solve using a uniform argument. Many of the about 20 surfaces that can be obtained this way were known to Karcher and Schoen, but some appear to be new.

Our method can also be used to describe more complicated surfaces of large genus where the period problem becomes high dimensional. Here, there are currently no elementary methods available to prove the existence of these surfaces. However, our method can still be easily implemented to yield numerical results and images.

\section{An equivariant Schwarz-Christoffel formula}

In this section, we will prove an equivariant version of the classical Schwarz-Christoffel formula for periodic polygons. The proof is a straightforward generalization of the classical case.

Instead of the monomial factors $z-p_i$ of the classical Schwarz-Christoffel formula, we will use $\vartheta$-factors $\vartheta(z-p_i)$ where
\[
\vartheta(z)
=\vartheta(z,\tau)
%=-\sum_{n=-\infty}^\infty e^{\pi i (n+\frac12)^2\tau+2\pi i (n+\frac12) (z+\frac12)}
=\sum_{n=-\infty}^\infty e^{\pi i (n+\frac12)^2\tau+2\pi i (n+\frac12) (z-\frac12)}
\]
is one of the classical Jacobi $\vartheta$-functions (\cite{mum1}). It is an entire function  with simple zeroes at the lattice points of the integer lattice spanned by $1$ and $\tau$.
It enjoys the following symmetries:

\begin{align*}
\vartheta(-z)&= -\vartheta(z),
\\
\vartheta(z+1)&=-\vartheta(z),
\\
 \vartheta(z+\tau)&=-e^{-\pi i \tau-2\pi i z} \vartheta(z).
\end{align*}

Further properties of $\vartheta(z)$ which we will need are that $\vartheta'(0)\ne 0$, and that $\vartheta(\bar z)=\overline{\vartheta(z)}$ for purely imaginary $\tau$, in particular $\vartheta(z)$ is real for real $z$. 
All these properties characterize $\vartheta(z)$ uniquely. Now we are ready to describe the image domains of our Schwarz-Christoffel formula:

\begin{definition}
\label{df:poly-arc}
A \emph{polygonal arc} is a  piecewise linear curve with a discrete vertex set.
%\authornote{In Definition \ref{df:poly-arc}, does ``simple'' mean no self-intersection?  
%            If so, then the last sentence in Definition \ref{df:ppp} does not make sense.  
%            If not, the angle condition \eqref{eq:angle1} can be $2k\pi$.
%           }
\end{definition}

\begin{definition}
\label{df:ppp}
A \emph{periodic polygon} $P$ is a simply connected domain in the plane with the following properties:
\begin{enumerate}
\item $P$ is bounded by two infinite polygonal arcs.
\item $P$  is invariant under a euclidean translation $V(z)=z+v$ for some $v\ne 0$.
\item The quotient $P/\langle V\rangle$ is conformally an annulus.
\end{enumerate}
\end{definition}

More generally, we will also allow as periodic polygons simply connected Riemann surfaces together with a flat structure
and two periodic polygons as boundary that are invariant under a holomorphic transformation preserving the
flat structure. This way periodic polygons can have vertices with interior angles larger than $2\pi$. 

We denote the vertices of the two polygonal arcs by $P_i$ and $Q_j$ so that $V(P_i) = P_{i+m}$ and $V(Q_j) = Q_{j+n}$ for all $i,j$ and fixed integers $m,n$. Denote the interior angles of the polygon at $P_i$ (resp. $Q_j$) by $\alpha_i$ (resp. $\beta_j$).
By assumption, these numbers are also periodic with respect to $m$ and $n$, respectively.

\begin{figure}[H] 
  \centering
  \includegraphics[width=4in]{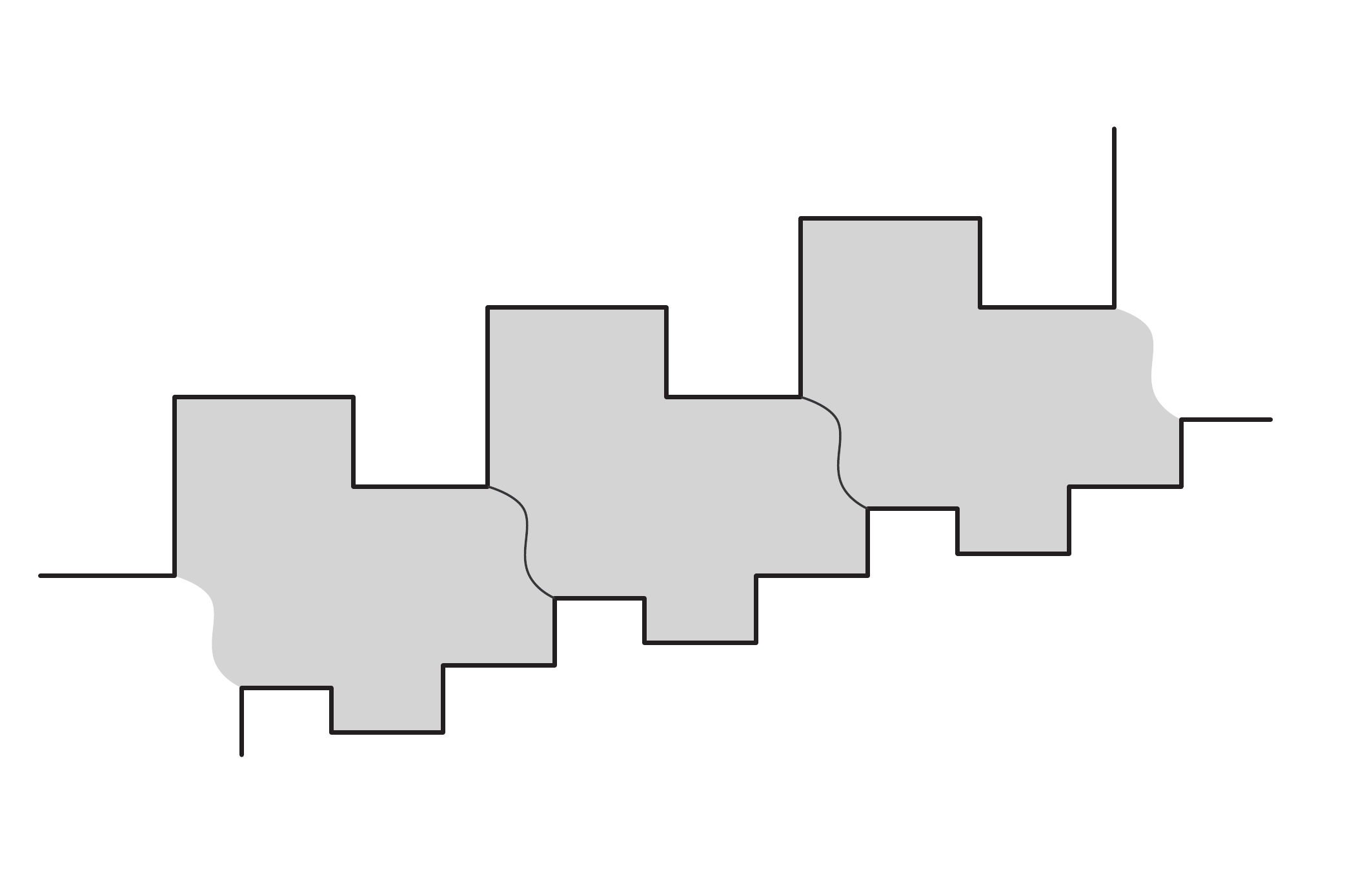} 
  \caption{A periodic polygon}
  \label{fig:frieze}
\end{figure}

Observe that 
\begin{equation}
\sum_{i=1}^m (\pi-\alpha_i)=0=\sum_{j=1}^n (\pi-\beta_j)
\label{eq:angle1}
\end{equation}
because both boundary arcs are invariant under a translation.

Let $P$ be a periodic polygon invariant under $V$. Let $d/2$ be the modulus of the annulus $P/\langle V\rangle$, and define  the strip $Z=\{z\in\C: 0<\im z<d/2\}$. The choice of $d$ makes 
the annuli  $Z/\langle z\mapsto z+1\rangle$ and  $P/\langle V\rangle$ conformally equivalent, and $d>0$  is uniquely determined this way. 
Moreover, we obtain a biholomorphic
map $f:Z\to P$ which is equivariant with respect to both translations:
\[
f(z+1)= V(f(z)).
\]
This map extends to a homeomorphism between the closures of $Z$ and $P$. In comparison to the classical Schwarz-Christoffel formula, $Z$ will play the role of the upper half  plane.

\begin{proposition}
Let $P$ be a periodic polygon, and $Z$ be the associated parallel strip as above. Then, up to scaling, rotating, and translating,
\[
f(z)=\int^z \prod_{i=1}^{m} {\vartheta(z-p_i)}^{a_i} \cdot \prod_{j=1}^{n} {\vartheta(z-q_j)}^{b_j}  
\]
is the biholomorphic map from $Z$ to $P$. Vice versa, for any choices of $p_i\in\R$ $(1\le i \le m)$, $q_j \in\R+di/2$ $(1\le j \le n)$ and $-1<a_i,b_j<1$ satisfying the angle condition
\begin{equation}
\sum_{i=1}^m a_i = 0 = \sum_{j=1}^n b_j, 
\label{eq:angle2}
\end{equation}
$f$ maps $Z$ to a (possibly immersed) periodic polygon.
\end{proposition}

\begin{proof}
Denote the preimages of $P_i$ and $Q_j$ under $f$ by $p_i$ and $q_j$. Let $\tau =i d$.
By the Schwarz reflection principle, $f$ can be holomorphically continued along any path in $\C$ avoiding the points $p_i +\Z \tau$ and    $q_j +\Z \tau$. It thus becomes a multivalued function
in the punctured plane $C_0 = \C - \{p_i +\Z \tau\}-\{q_j +\Z \tau\}$. As two consecutive reflections of the polygon $P$ result in a euclidean motion, the pre-Schwarzian derivative  $\eta=(f''/f')\, dz = d\, \log(f'(z)) $ of $f$ is a well-defined
holomorphic 1-form in $C_0$. This 1-form $\eta$ is not only invariant under the translation by $1$ (because of the equivariance of $f$) but also under translation by $\tau$ (by the definition of the pre-Schwarzian).

It thus descends to a holomorphic 1-form on the quotient $ C_0/\langle 1,\tau \rangle$. This quotient surface is conformally a rectangular torus punctured at finitely many points.

Moreover, $\eta$ extends meromorphically to a 1-form (also denoted by $\eta$)  on the rectangular torus $T= \C/\langle 1,\tau \rangle$
with first order poles at $p_i$ and $q_j$.  By the same argument as in the classical Schwarz-Christoffel formula, the residues at $p_i$ (resp. $q_j$) are  $a_i=\alpha_i/\pi-1$ (resp. $b_j=\beta_j/\pi-1$).

Using $\vartheta(z)=\vartheta(z,\tau)$, it is easy to write down such a 1-form as
\[
\eta_0 = \sum_{i=1}^{m} a_i \frac{\vartheta'(z-p_i)}{\vartheta(z-p_i)} \,dz+ \sum_{j=1}^{n} b_j \frac{\vartheta'(z-q_j)}{\vartheta(z-q_j)} \, dz.
\]
We claim that in fact $\eta=\eta_0$. 

First note that it follows from the transformation laws of $\vartheta$ and the angle condition \eqref{eq:angle1} or \eqref{eq:angle2} that $\eta_0$ is indeed elliptic, because
with $h(z)=\vartheta'(z)/\vartheta(z)$ we have
\[ 
  h(z+1) =  h(z) , \qquad
  h(z+\tau) = -2\pi i + h(z).
\]

Thus $\eta-\eta_0$ is a holomorphic 1-form on  $T$, i.e. a multiple of $dz$.
To see that $\eta=\eta_0$, we integrate both over the cycle $\gamma$ on $T$ homologous to $\tau /2+[0,1]$:

As $f(z+1)=f(z)+v$ we have $f'(z+1)=f'(z)$ and thus $\int_{\gamma } \eta=0$.

On the other hand, as $h(z)$ is well-defined on the cylinder $\C/\langle z\mapsto z+1\rangle$ (with generating cycle $\gamma$),
%\[
%\int_{\gamma} h(z)\, dz= \log (-1)
%\]
%This implies
\[
\int_{\gamma } h(z-p_i)\, dz=  \int_{\gamma } h(z-q_j)\, dz
\]
and thus the angle condition \eqref{eq:angle2} implies that 
$\int_{\gamma } \eta_0=0$. Therefore, the holomorphic 1-form $\eta-\eta_0$ has $\gamma $-period 0, and must vanish identically.

By integrating $\eta_0$ we get
\[
\log(f'(z))=\int^z \eta_0 = \sum_{i=1}^{m} a_i\log({\vartheta(z-p_i)}) + \sum_{j=1}^{n} b_j \log({\vartheta(z-q_j)})  +C
\]
or
\[
f'(z)=e^C \prod_{i=1}^{m} {\vartheta(z-p_i)}^{a_i} \cdot \prod_{j=1}^{n} {\vartheta(z-q_j)}^{b_j} ,
\]
which proves our formula.

The other direction follows by the Schwarz reflection principle and the local analysis at the $p_i$ and $q_j$ as in the proof of the standard Schwarz-Christoffel formula.
\end{proof}

\begin{remark}
The same proof works also for infinite polygons that are invariant under a rotation (possibly by an irrational angle). These polygons will only be immersed, but this doesn't cause a problem. The only change required is that the angle condition \eqref{eq:angle1} relaxes to
\[
\sum_{i=1}^m (\pi-\alpha_i)=\sum_{j=1}^n (\pi-\beta_j), 
\]
which is all that was needed in the above proof.

With some more notational effort, it is also possible to discuss more general  polygons where the interior angles do not necessarily lie between $0$ and $2\pi$.

However, we will not need such polygons in this paper.
\end{remark}

\section{Minimal surfaces defined on parallel strips}

In this section we discuss simply connected, periodic minimal surfaces with two boundary components whose Weierstrass data can be given by the integrands of our Schwarz-Christoffel formula. These surfaces constitute  the building blocks for the triply periodic surfaces constructed in the following sections.

Any conformally parametrized minimal surface $\Sigma$ in euclidean space can locally be given by the Weierstrass representation
\[
z\mapsto \re \int^z\left( \frac{1}{2}\left(\frac{1}{G}-G\right),\,
                \frac{i}{2}\left(\frac{1}{G}+G\right),\,
                1 \right) \, dh,
\]
where $G$ is a meromorphic function and $dh$ is a holomorphic 1-form.  
$G$ can be identified with the Gauss map of $\Sigma$ via the stereographic projection, and $dh$ is called the height differential of $\Sigma$.  
The triple $(\Sigma,G, dh)$  are the Weierstrass data of the minimal surface.

For a positive real number  $d$, let $Z=\{z\in\C: 0<\im z< d/2\}$ a strip domain. Introduce
$\tau = d \cdot i$ and the rectangular torus $T= \C/\langle 1,\tau \rangle$.

For $m,n,\ge 0$ consider points $p_i\in (0,1), i=1,\ldots,m$ and $q_j\in (0,1)+\tau /2, j=1,\ldots,n$. Extend these to periodic sets of points by imposing the condition
$p_{i+m}=p_i+1$, $q_{j+n}=q_j+1$ for all $i,j$.

Consider the minimal surface $\Sigma$ defined on  $Z$ by the Weierstrass data
\[
G(z)= \prod_{i=1}^{m} {\vartheta(z-p_i)}^{a_i} \cdot \prod_{j=1}^{n} {\vartheta(z-q_j)}^{b_j}  
\]
and $dh =dz$ on $\Sigma = T$.

For the exponents we assume that $-1<a_i,b_j<1$ for all $i,j$ --- this ensures that all interior angles of  the Schwarz-Christoffel polygons are between $0$ and $2\pi$.
We also assume the angle condition \eqref{eq:angle2}. Then we have:

\begin{proposition}
\label{prop:corner} 
$\Sigma$ is a simply connected minimal surface with two boundary components lying in  a finite number of vertical symmetry planes. These planes meet at angles $\pi a_i$  at the image of $p_i$ and $\pi b_j$  at the image of $q_j$. Furthermore, the surface is invariant   under the vertical translation $x_3\mapsto x_3+1$.
\end{proposition}

\begin{proof}
Define the following functions on $Z$:
\[ 
 \Phi_1(z) =  \int^z G \,dh , \qquad 
 \Phi_2(z) =  \int^z \frac1{G} \,dh .
\]

By our Schwarz-Christoffel formula, both $\Phi_1$ and $\Phi_2$ map $Z$ to periodic polygons. We first claim that corresponding oriented segments of these polygons represent conjugate directions (they are most likely of different lengths, though). This statement remains true or false if we multiply $G$ by a fixed factor $e^{i \phi}$.

Thus we can assume without loss of generality that the segment under consideration of $Z$ is mapped by $\Phi_1$ to a positively oriented horizontal segment. This means that $G$ must be positive and real on that segment. Consequently, $\Phi_2$ maps the same segment also to a positively oriented horizontal segment. This implies our claim.

Now let 

\[
F(z)= \re \int^z \frac{1}{2}\left(\frac{1}{G}-G\right) \, dh 
   +i\re \int^z \frac{i}{2}\left(\frac{1}{G}+G\right) \, dh .
\]

This is the orthogonal projection of $\Sigma$ to the $x_1x_2$-plane.
Then
\[
2F(z) = \overline{\Phi_2(z)} - \Phi_1(z).
\]
We normalize the integration constants for $\Phi_1$ and $\Phi_2$ and rotate $G$ if necessary so that the  image segment under consideration lies  on the real axis. Then it follows that
\[
F(\bar z) =\overline{F(z)}
\]
so that the image curve of the segment under consideration is indeed a planar symmetry curve.

To compute the angle between symmetry planes at the corners $p_i$ (or $q_j$), note that these corners are points where the Gauss map becomes vertical. Furthermore, just before and after $p_i$  (or $q_j$) the Gauss map moves along straight lines in $\C$ through $0$. By the explicit formula for $G(z)$, it changes by the factor $(-1)^{a_i}$ when crossing $p_i$ or  $(-1)^{b_j}$ when crossing $q_j$. This implies the claimed angles.

The invariance under the vertical translation is an immediate consequence of the invariance of $\Phi_1$ and $\Phi_2$ under the translation $z\mapsto z+1$ and the definition $dh=dz$.

\end{proof}

\begin{definition} For an interval $[p_i,p_{i+1}]$ or $[q_j,q_{j+1}]$, we denote the symmetry plane in which the corresponding planar symmetry curve lies by $\Pi_{[p_i,p_{i+1}]}$ or $\Pi_{[q_j,q_{j+1}]}$. We orient these plane by insisting that its normal vector points away from the minimal surface.
\end{definition}

The following corollary follows from standard  properties of conjugate minimal surfaces.
\begin{corollary}
The conjugate surface of $\Sigma$ is a simply connected minimal surface bounded by two spatial polygonal arcs. The angles at the images of $p_i$ $($resp. $q_j)$ are $\pi a_i$ $($resp. $\pi b_j)$.
\end{corollary}

\section{The Basic Examples}
\label{sec:basic}

\def\fw{1.5in}
\begin{figure}[H]
 \begin{center}
   \subfigure[(2,4,4)]{\label{1(2,4,4)}\includegraphics[width=\fw]{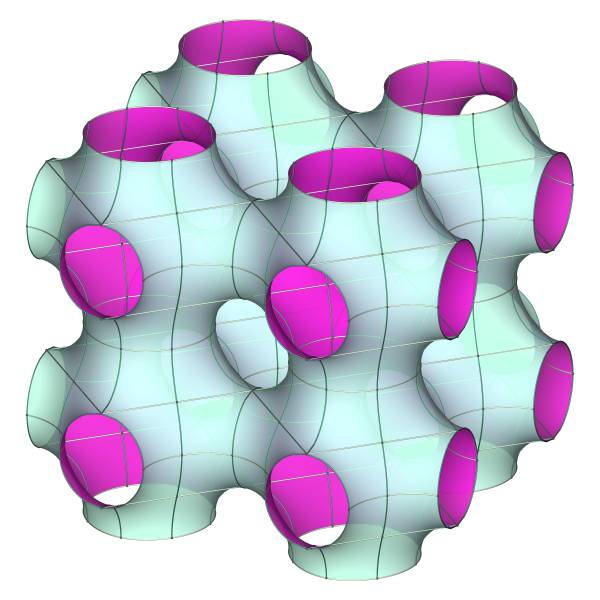}}
   \subfigure[(3,3,3)]{\label{1(3,3,3)}\includegraphics[width=\fw]{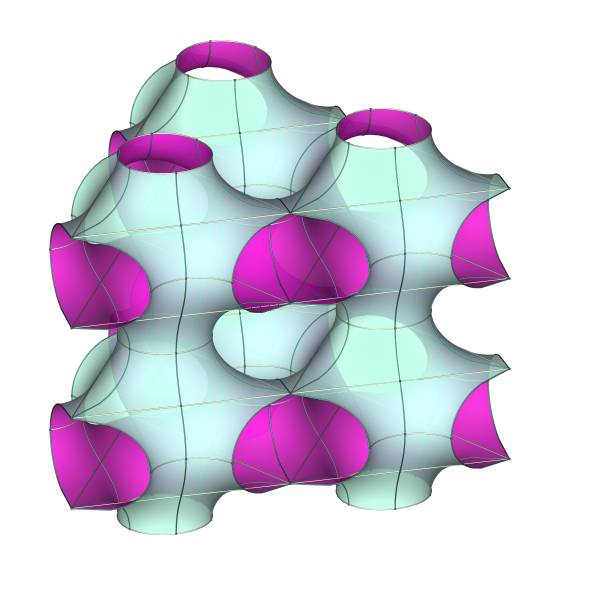}}
   \subfigure[(4,2,4)]{\label{1(4,2,4)}\includegraphics[width=\fw]{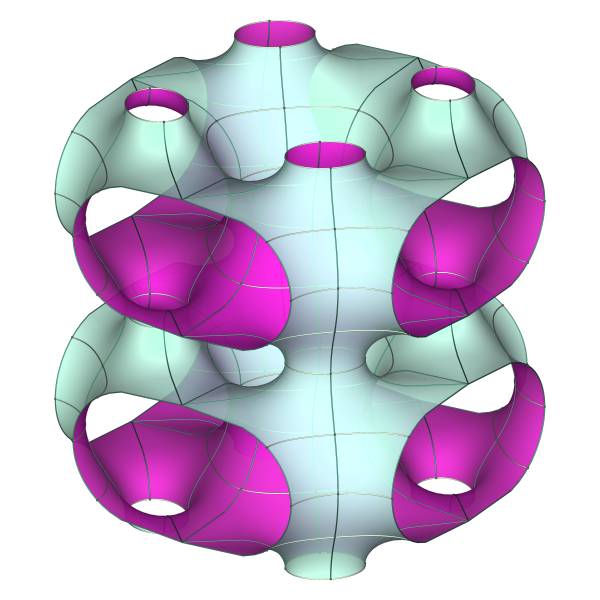}}\\
   \subfigure[(4,4,2)]{\label{1(4,4,2)}\includegraphics[width=\fw]{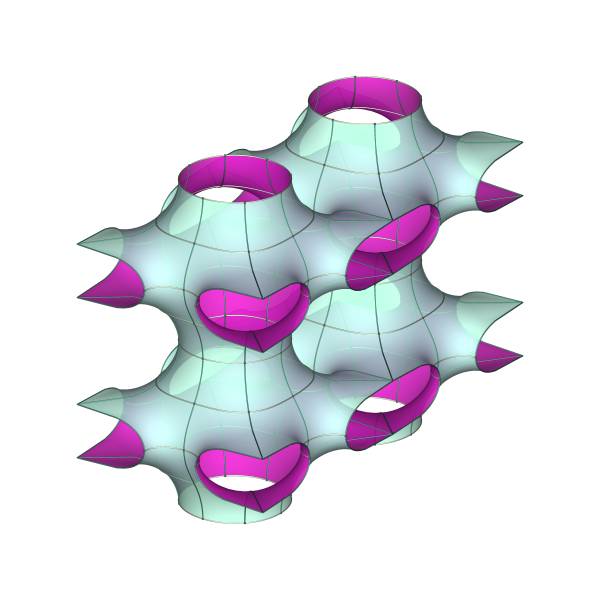}}
   \subfigure[(2,3,6)]{\label{1(2,3,6)}\includegraphics[width=\fw]{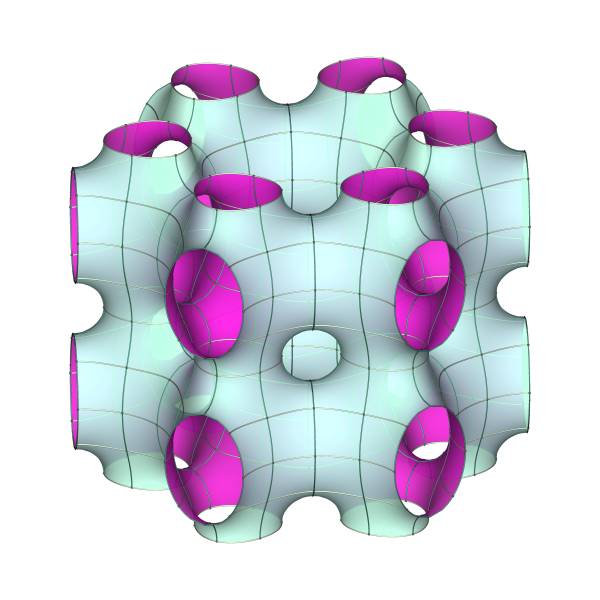}}
   \subfigure[(2,6,3)]{\label{1(2,6,3)}\includegraphics[width=\fw]{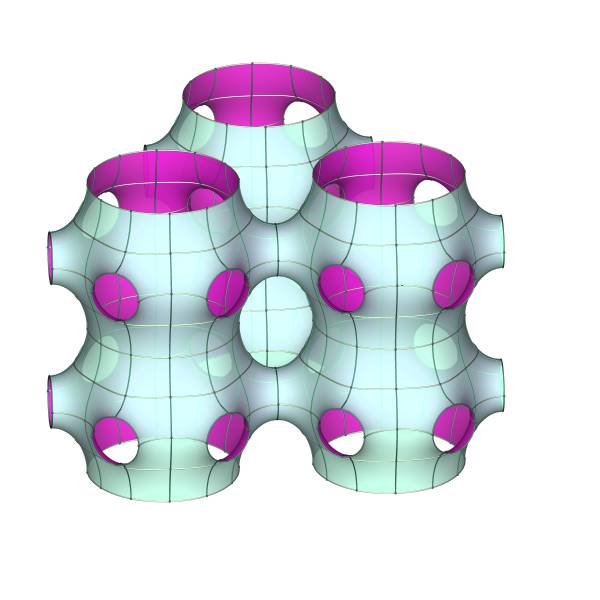}}\\
   \subfigure[(3,2,6)]{\label{1(3,2,6)}\includegraphics[width=\fw]{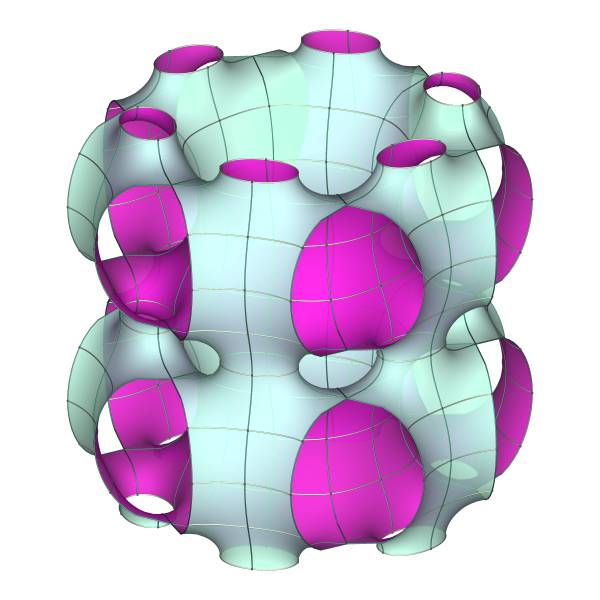}}
   \subfigure[(3,6,2)]{\label{1(3,6,2)}\includegraphics[width=\fw]{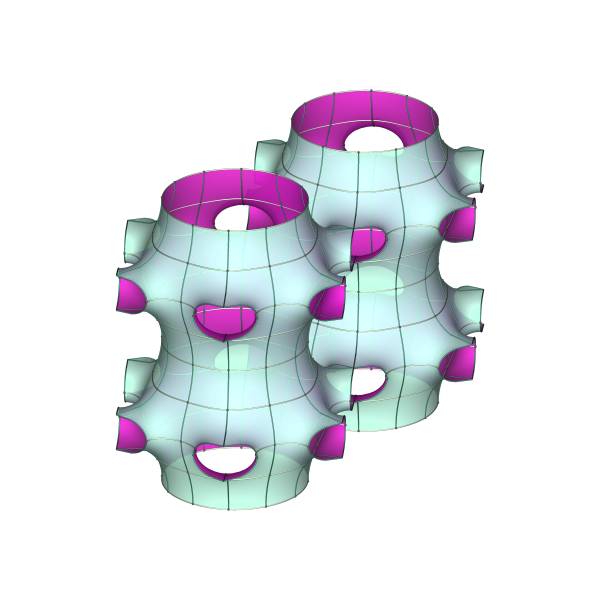}}
   \subfigure[(6,2,3)]{\label{1(6,2,3)}\includegraphics[width=\fw]{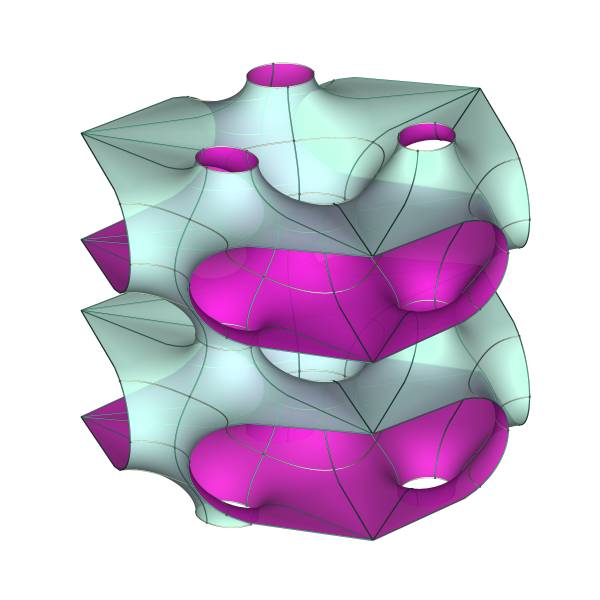}}\\
   \subfigure[(6,3,2)]{\label{1(6,3,2)}\includegraphics[width=\fw]{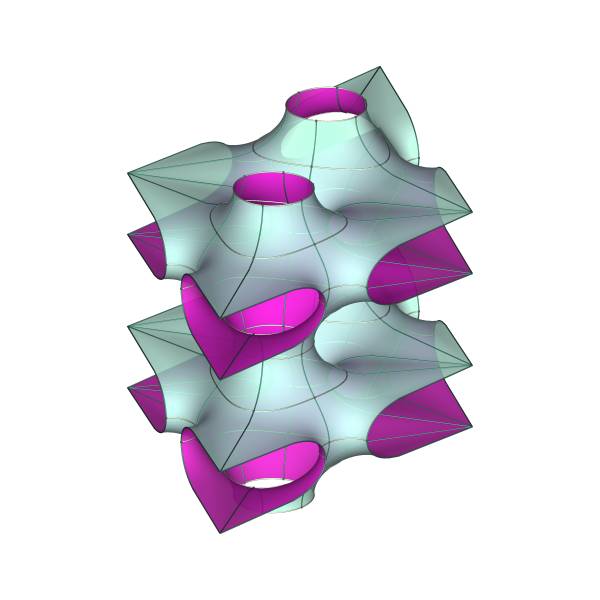}}
 \end{center}
 \caption{The Basic Examples}
 \label{fig:basicexamples}
\end{figure}

The simplest example of a periodic polygon has just two vertices along one edge and none along the other. We will now discuss the resulting minimal surface $\Sigma$. By translating $Z$, we can assume that $p_1=-p$ and $p_2=+p$. 
Furthermore, by the angle condition \eqref{eq:angle2}, $a_1=a=-a_2$. Thus the Gauss map is given by
\[
G(z)=\left(
\frac{\vartheta(z-p)}{\vartheta(z+p)}\right)^{a} 
\]
for some $0<a<1$.

With this normalization, the minimal surface $\Sigma$ becomes symmetric with respect to a reflection at the imaginary axis in the domain which corresponds to a reflection at a horizontal plane in space, which we can assume to be the $x_1x_2$-plane.

By proposition \ref{prop:corner},  $\Sigma$ has two boundary arcs. The one corresponding to $\im z = \im\tau /2$ lies in a single coordinate plane $\Pi$, while the other switches between coordinate planes $\Pi_{[-p,p]}$ and $\Pi_{[p,-p]}$, making an angle $a\pi$. Our next goal is to compute the angle between $\Pi$ and $\Pi_{[-p,p]}$:

\begin{proposition}
\label{prop:bottom}
The angle $\alpha_0$ between $\Pi_{[-p,+p]}$ and $\Pi$ is equal to $\alpha_0=\pi a(2p-1)$.
\end{proposition}
\begin{proof}
By the quasiperiodicity of $\vartheta$, 
\begin{align*}
\vartheta(p+\tau/2) ={}& \vartheta(p-\tau/2+\tau)
\\
={}&- \vartheta(p-\tau/2)e^{-\pi i \tau -2\pi i( p-\tau/2)}
\\
={}& \vartheta(-p+\tau/2)e^{ -2\pi i p}
\end{align*}

Note that along the segment $[0,\tau/2]$ the Gauss map is horizontal. Thus
\begin{align*}
\arg(G(\tau/2))-\arg(G(0)) 
={}&\arg\left\{ \left(
\frac{\vartheta(\tau/2-p)}{\vartheta(\tau/2+p)}\right)^{a} \left(
\frac{\vartheta(p)}{\vartheta(-p)}\right)^{a} \right\}
\\
={}& 2\pi  p a-\pi  a.
\end{align*}
This computation shows the last equality modulo $2\pi a$. As both sides  are continuous in $p$ for $0<p<1$, the claim follows because it is true when $p=1/2$ and $G(z)\equiv (-1)^a$.
\end{proof}

Our next goal is to obtain a triply periodic surface from $\Sigma$ by repeated reflections at the three symmetry planes. In order for the surface to be embedded, the group generated by these reflection must be a euclidean triangle group.

There are three such groups, denoted by $\Delta(2,3,6)$, $\Delta(2,4,4)$, and $\Delta(3,3,3)$ where $\Delta(r,s,t)$ corresponds to the group generated by reflecting at the edges of a triangle with angles $\pi/r$, $\pi/s$, and $\pi/t$. We assign a specific such triangle group to $\Sigma$ as follows:

The angle at $+p$ is to be $\pi/r$, the angle between $\Pi_{[-p,p]}$ and $\Pi$ equals $\pi/s$, and the angle between $\Pi_{[p,-p]}$ and $\Pi$ shall be $\pi/t$. 
By Proposition \ref{prop:corner}, the first condition forces $a= (r-1)/r$, while the second determines 
\[
p=\frac{-r - s + r s}{2 (r-1) s}
\]
by Proposition \ref{prop:bottom}. 

This allows for 10 possibilities. However, due to the reflectional symmetry at the $x_1x_2$-plane, any surface  corresponding to $\Delta(r,s,t)$ becomes one corresponding to $\Delta(r,t,s)$, by turning it upside down. Thus there are only 5 distinct cases. They all correspond to surfaces known to H.~Schwarz (\cite{schw1}) or A.~Schoen (\cite{sch1}).

The following table lists all cases with the naming convention of A.~Schoen, as well as the value of $p$.

%\begin{center}
%\begin{tabular}{l|cccc}
%Name          & $(r,s,t)$ & $p$ & $\alpha_1(a)=\pi/s$ & $\alpha_2(a)=-\pi/t$ \\ \hline
%Schwarz P     &  (2,4,4)  & 1/4 & $45$          & $-45$         \\
%Schoen H'-T   &  (2,6,3)  & 1/3 & $30$          & $-60$         \\
%Schoen H'-T   &  (2,3,6)  & 1/6 & $60$          & $-30$         \\ \hline
%Schwarz H     &  (3,3,3)  & 1/4 & $60$          & $-60$         \\
%Schoen H''-R  &  (3,2,6)  & 1/8 & $90$          & $-30$         \\
%Schoen H''-R  &  (3,6,2)  & 3/8 & $30$          & $-90$         \\ \hline
%Schoen S'-S'' &  (4,4,2)  & 1/3 &$45$          & $-90$         \\
%Schoen S'-S'' &  (4,2,4)  & 1/6 & $90$          & $ 45$         \\ \hline
%Schoen T'-R   &  (6,2,3)  & 1/5 & $90$          & $-60$         \\
%Schoen T'-R   &  (6,3,2)  & 3/10& $60$          & $-90$
%\end{tabular}
%\end{center}

\begin{center}
\begin{tabular}{l|cc}
Name          & $(r,s,t)$ & $p$  \\ \hline
Schwarz P     &  (2,4,4)  & 1/4          \\
Schoen H'-T   &  (2,6,3)  & 1/3         \\
Schoen H'-T   &  (2,3,6)  & 1/6          \\ \hline
Schwarz H     &  (3,3,3)  & 1/4          \\
Schoen H''-R  &  (3,2,6)  & 1/8         \\
Schoen H''-R  &  (3,6,2)  & 3/8          \\ \hline
Schoen S'-S'' &  (4,4,2)  & 1/3         \\
Schoen S'-S'' &  (4,2,4)  & 1/6          \\ \hline
Schoen T'-R   &  (6,2,3)  & 1/5          \\
Schoen T'-R   &  (6,3,2)  & 3/10
\end{tabular}
\end{center}

All these surfaces come in a 1-parameter family where $\tau\in i \R^+$ is the parameter.

The embeddedness of all surfaces is easiest seen by the conjugate surface method explained in \cite{ka5,ka6}. 

\section{The general symmetric case}

We will now discuss minimal surfaces related to periodic polygons with more corners, restricting our attention to the following symmetric case: We assume that the surface is symmetric with respect to the $x_1x_2$-plane. We can assume that this reflection is realized in the domain by a reflection at the imaginary axis.

This implies that the  Gauss map is symmetric with respect to the imaginary axis
as in the previous examples, i.e. both the $p_i\in\R$, ($1\le i \le m$) and the $q_j\in\R +\tau /2$, ($1\le j \le n$) are symmetric with respect to the imaginary axis, and the exponents have opposite signs. In other words, we assume that

\begin{equation}
G(-z)=\frac1{G(z)}.
\label{eqn:symmetry}
\end{equation}

Relabel the $p_i$ and $q_j$ so that $0<p_1<p_2<\cdots <p_{m'}<1/2$, ($m'=m/2$) and $0<\re(q_1)<\re(q_2)<\cdots<\re(q_{n'})<1/2$, ($n'=n/2$). This way, the interval $[-p_1,p_1]$ corresponds to an edge of the periodic polygons, and is itself symmetric with respect to the assumed additional reflectional symmetry.
Again, we need a formula for the angle between symmetry planes in the two different boundary components:

\begin{proposition}
\label{prop:bottom2}
The angle $\alpha_0$ between $\Pi_{[-p_1,p_1]}$ and $\Pi_{[-q_1,q_1]}$ is equal to 
\[
\alpha_0 =\pi \left(\sum_{i=1}^{m'} a_i(2p_i-1)+\sum_{j=1}^{n'} b_j(2\re(q_j)-1)\right)\ .
\]

\end{proposition}
\begin{proof}
In this case  the Gauss map is 
\[
G(z)=
\prod_{i=1}^{m'}
\left(\frac{\vartheta(z-p_i)}{\vartheta(z+p_i)}\right)^{a_i}
\prod_{j=1}^{n'}
\left(\frac{\vartheta(z-q_j)}{\vartheta(z+q_j)}\right)^{b_j}\quad,
\]
the factors of which were dealt with in the proof of proposition \ref{prop:bottom}.
\end{proof}

The additional symmetry \eqref{eqn:symmetry} about the $x_1x_2$-plane implies an important symmetry property of the periodic polygons given by  $G dh$ and $1/G dh$:

\begin{proposition}
\label{prop:symmetry}
Suppose that $G(z)$ is symmetric as above. Then the periodic polygons $\Phi_1(Z)$ and $\Phi_2(Z)$ are symmetric to each other via a reflection about a vertical line.
\end{proposition}

\begin{proof} This follows from the computation
\begin{align*}
\Phi_1(-\bar z) ={}& \int_0^{-\bar z} G(w) \, dw 
= \overline{\int_0^{- z} {G(w) \, dw }}
= -\overline{\int_0^{ z} {G(-w) \, dw }} \\
={}& -\overline{\int_0^{ z} {\frac1{G(w)} \, dw }}
= -\overline{\Phi_2(z)}\ .
\end{align*}
\end{proof}

\section{Four corners, all in one boundary component}

The next complicated case after the basic case (section \ref{sec:basic}) with two corners allows for four corners. These can either lie in a single boundary component, or be divided into two corners for each component.
{In this section, we will discuss four corners lie in a single boundary component.} 

To simplify the notation, we label the $p_i$ as
\[
(p_1,p_2,p_3,p_4)= (-q,-p,p,q)
\]
and the exponents as
\[
(a_1,a_2,a_3,a_4)=(-b,-a,a,b).
\]

By taking an (unbranched) double cover over our basic examples by doubling the vertical period, one can always obtain surfaces of this type. In this case, $q={1-p}$ and $b=-a$. Observe that  the exponents of consecutive points alternate in sign $(+,-,+,-)$.
We do not have any hope of finding other solutions in this symmetry class with the same sign pattern of the exponents.

\begin{figure}[H]
\centering
  \includegraphics[width=3in]{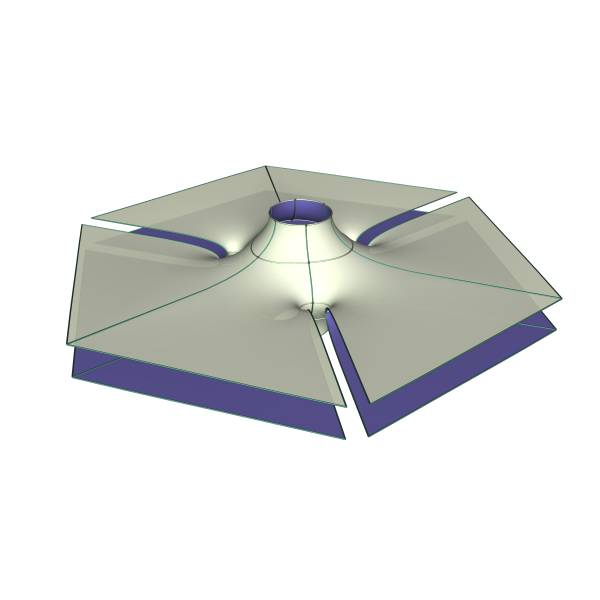} 
  \caption{A $(3,3,3)$-candidate with unclosable period}
  \label{fig:unclosed}
\end{figure}

A second type of plausible candidates arises by letting the exponents have signs $(-,-,+,+)$. We will now show that in this case, the period problem can never be solved. The picture  in figure \ref{fig:unclosed} illustrates that the period gap can be made arbitrarily small by letting $\tau\to0$.

To prove the  impossibility of a solution to the period problem, first note that the image of the lower edge of $Z$ determines the shape of the reflection triangle completely. This makes it necessary that the  upper edge of $Z$  is mapped  into a plane over one of these edges, and we can assume without loss of generality it to be the image edge of $[-p,p]$. In particular, these image edges become parallel, and the angle condition becomes
\[
2(ap+bq)=a+b-1 \ .
\]
This given, the  period condition requires the two image edges to lie above each other,  or
\[
\re \int_0^{\tau/2}{\frac{i}{2}\left(\frac{1}{G}+G\right)\,dh} = 0.
\]
In terms of the periodic polygon $\Phi_1(Z)$, this is equivalent to the condition that the two segments are collinear.

\begin{figure}[H]
\centering
  \includegraphics[width=5in]{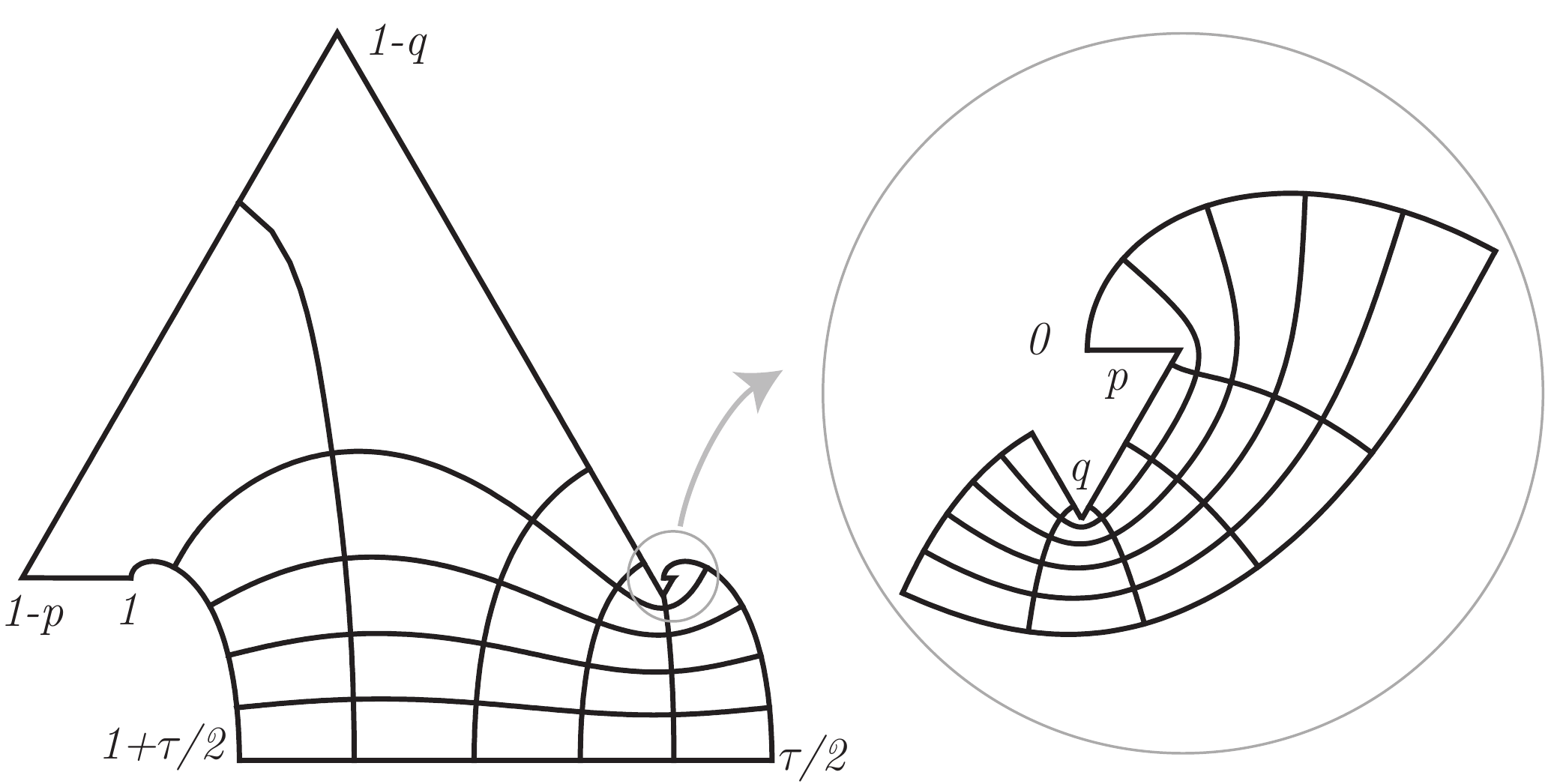} 
  \caption{Emerging branched point obstructs the solution of the period problem}
  \label{fig:branched}
\end{figure}

Figure \ref{fig:branched} shows the image of $[0,1]\times [0,\tau/2]$ under $\Phi_1$. We have labeled the image vertices by using the points in $Z$ for simplicity. The image to the left shows the entire image, the image to the {right} is a zoomed in portion as indicated. Without loss of generality, we can assume that $[0,p]$ and thus also $[1+\tau/2,\tau/2]$ are horizontal.

By assumption, the interior angles at $0$ and $\tau/2$ are ${\pi /2}$, while the angles at $p$ and $q$ are both larger than ${\pi}$, as $a,b>0$. This forces $q$ to lie below $[0,p]$. On the other hand, the period condition forces $[0,p]$ to be collinear with $[1+\tau/2,\tau/2]$. This, however, is a contradiction, as there are no interior branched points of $\Phi_1$ in $Z$.

\section{Two corners in each boundary component}

In this section, we will discuss symmetric polygons where each boundary component has just two corners, i.e. $m=n=2$. In contrast to the previous sections, we will obtain new examples.

Without loss of generality we can assume that $p_1=p=-p_2\in(0,1/2)$ and 
$q_1=q=-q_2\in(0,1/2)+\tau/2$.

There are two qualitatively distinct cases, depending on whether $a=a_1$ and $b=b_1$ have the same or opposite signs.

\subsection{Exponents have equal signs}

Without loss of generality, we can assume that $a,b<0$.

In this case, projecting the two boundary arcs into the $x_1x_2$-plane gives two ``hinges'' with angles $\pi a$ and $\pi b$. 
We denote the hinge  containing the image of $p$ by $H_p$, and the hinge containing $q$ by $H_q$.

We want these hinges to be part of our reflection group triangles. The requires $a$ and $b$ to be of the form
  $ a = -(r-1)/r$ and $b = -(s-1)/s$ with $r,s\in\{2,3,4,6\}$. 

\def\fw{2in}
\begin{figure}[H]
 \begin{center}
   \subfigure{\includegraphics[width=\fw]{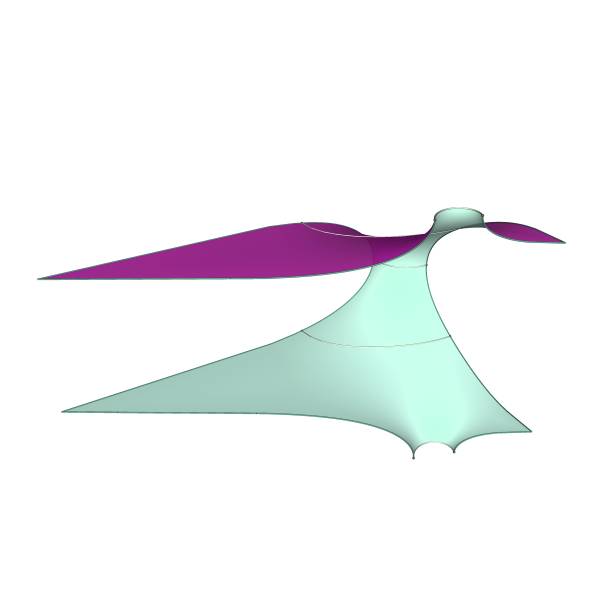}}
   \subfigure{\includegraphics[width=\fw]{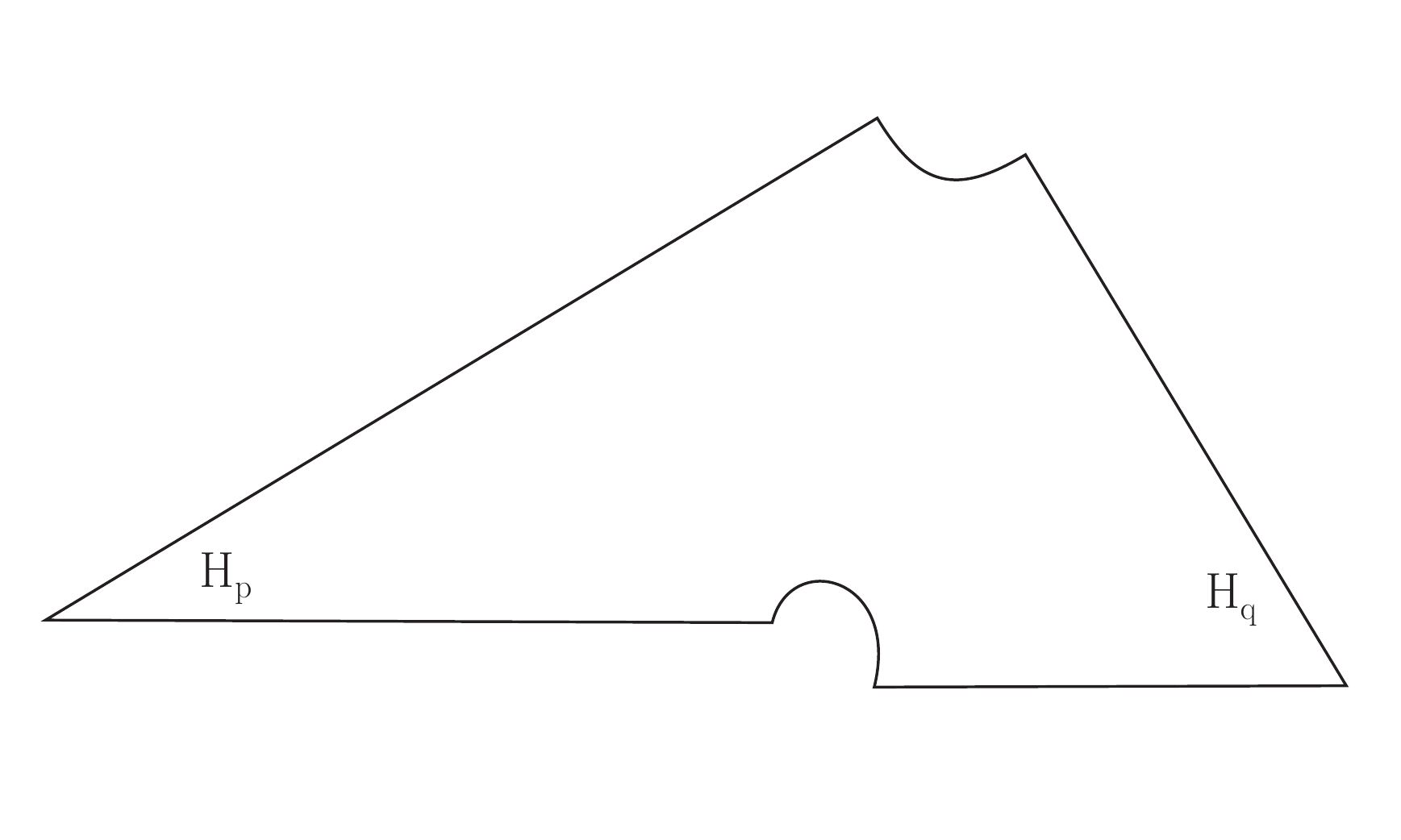}}
     \end{center}
 \caption{Minimal $(3,6,2)$-patch with schematic hinges of the symmetry planes }
 \label{fig:hinges1}
\end{figure}

As the hinges have four edges all together, two of these edges must lie on one side of a reflection group triangle. By relabeling the vertices, if necessary, we can thus assume that the planes $\Pi_{[-p,p]}$ and $\Pi_{[-q,q]}$ are parallel.

But proposition \ref{prop:bottom2}, this means that 
\begin{equation}
\label{eqn:parallel}
a(2p-1)+ b(2 \re(q)-1)=1 .
\end{equation}

This constraint between $p$ and $q$ guarantees that the planes $\Pi_{[-p,p]}$ and $\Pi_{[-q,q]}$ are parallel, but we need them to be equal. We will show now that we can always adjust $p$ to make this happen:

\begin{theorem}
Given any $\tau\in i  \R^+$ and $-1<a,b<0$ , there are $0<p,\re(q)<1/2$ so that the two hinges $H_p$ and $H_q$  line up as a triangle. 
\end{theorem}

\begin{proof}
To show this, we have to adjust $p$ as to
satisfy period condition:

\[
\int_0^{\tau} G \, dh = \overline{\int_0^{\tau} \frac1G \, dh}.
\]

Because of the horizontal planar symmetry and proposition \ref{prop:symmetry}, this is equivalent to 
\[
\int_0^{\tau} G \, dh = 0.
\]

For the periodic polygon $\Phi_1(Z)$ this means that the image edges $\Phi_1([-p,p])$ and $\Phi_1([-q,q])$ need to be collinear.

Without loss of generality, we can assume that these edges are horizontal. By \eqref{eqn:parallel}, they are already parallel.

We will now show that for $p$ close to $0$ (so that, by equation \eqref{eqn:parallel}, $\re(q)$ is close to $(a+b+1)/2b$), $\Phi_1([-p,p])$ is below $\Phi_1([-q,q])$, while for $\re(q)$ close to $0$ (so that $p$ is close to $(a+b+1)/2a$), $\Phi_1([-p,p])$ is above $\Phi_1([-q,q])$. Then the claim follows from the intermediate value theorem.

\def\fw{1.6in}
\begin{figure}[H]
 \begin{center}
 \subfigure[$p\to 0$]{\label{interp}\includegraphics[width=\fw]{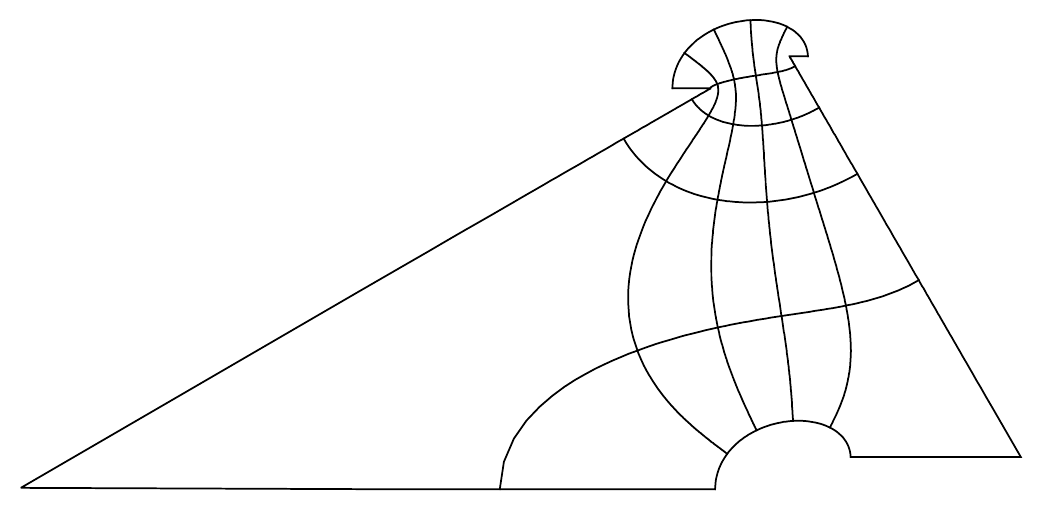}}
  \subfigure[period problem solved]{\label{intersol}\includegraphics[width=\fw]{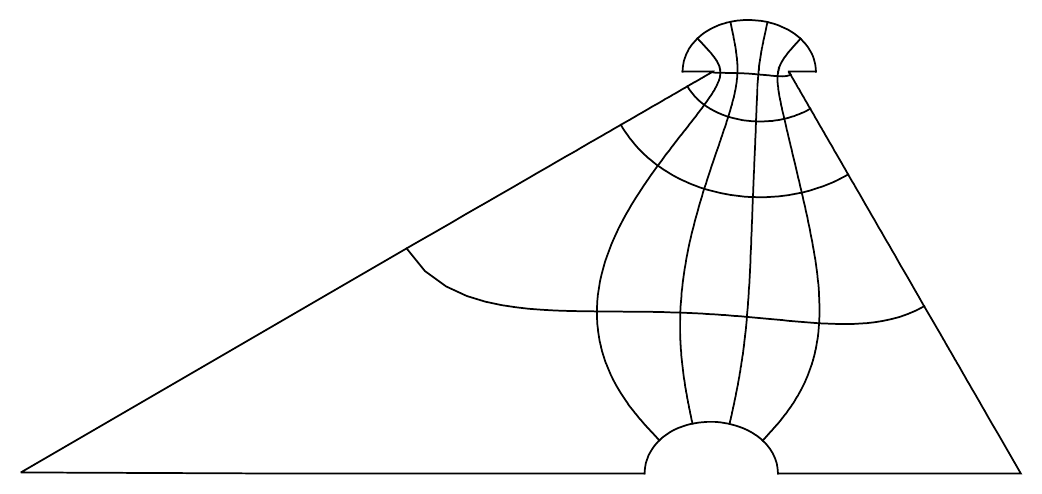}}
   \subfigure[$q\to \tau /2$]{\label{interq}\includegraphics[width=\fw]{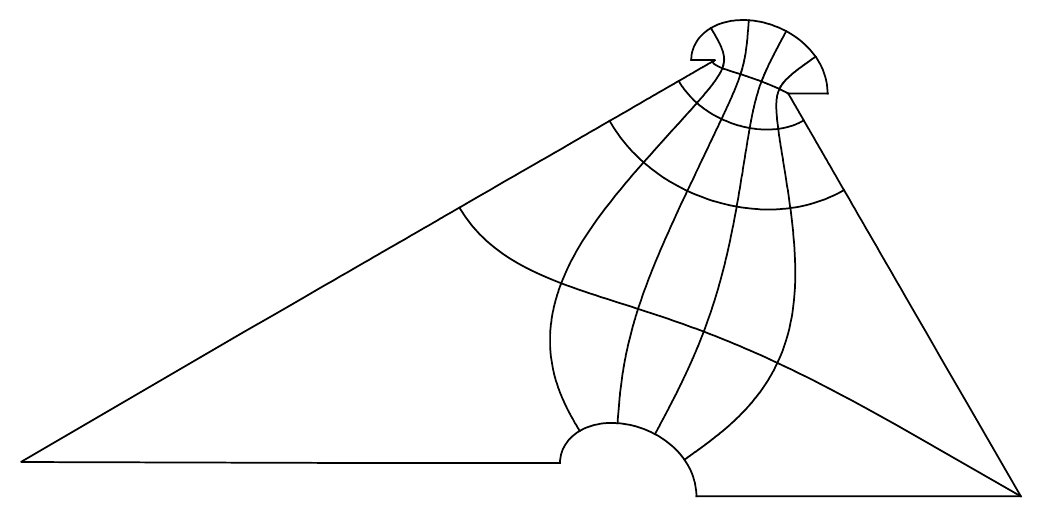}}
 \end{center}
 \caption{Intermediate value argument}
 \label{fig:intermediate}
\end{figure}

Let's first consider the case $q\to \tau /2$. To show that $\Phi_1([0,p])$ is below $\Phi_1([\tau /2,q])$, we consider the limit case.

\begin{figure}[H] %  figure placement: here, top, bottom, or page
  \centering
  \includegraphics[width=3in]{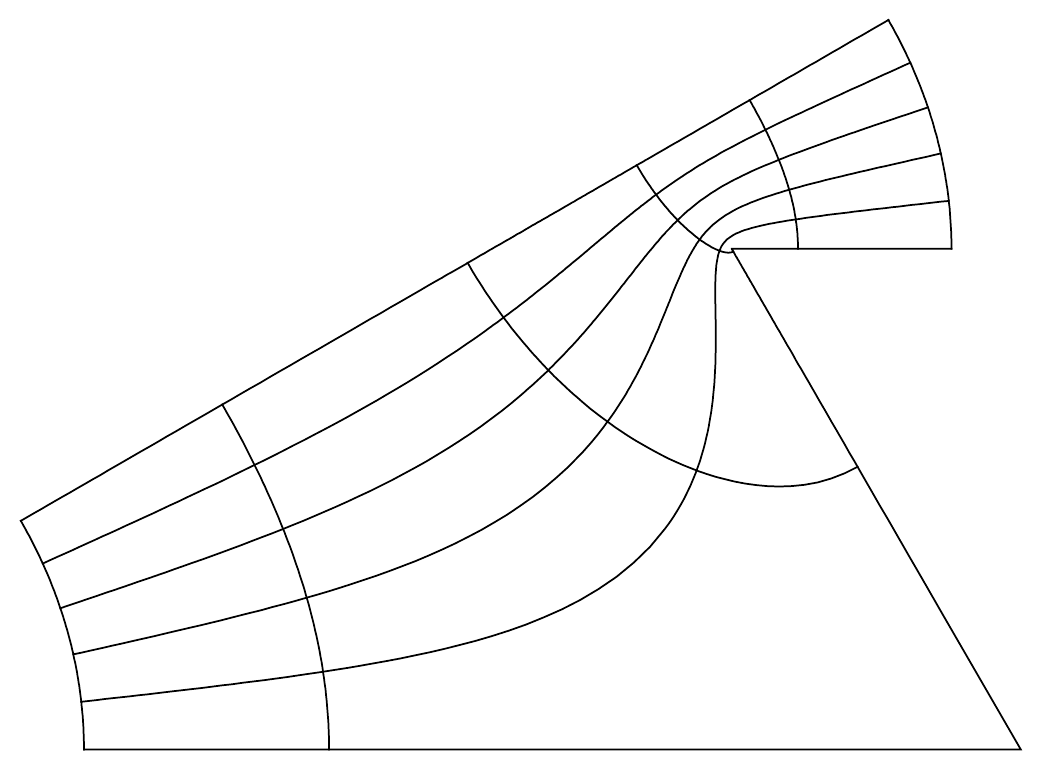} 
  \caption{Limit of the domains for $q\to \tau /2$}
  \label{fig:LimitDomain}
\end{figure}

As period integrals are continuous in the parameters, it is sufficient to show that for 
\[
\Phi_1(z) = \int^z 
\left(
\frac{\vartheta(w-p)}{\vartheta(w+p)}
\right)^{a} \,dw
\]
satisfies that $\Phi_1(\tau /2)$ lies above $\Phi_1(0)$. We will show that the height of the curve $\Phi_1(t\in [0,\tau /2])$ is strictly increasing with the parameter $t$. By conformality, this curve is perpendicular to the straight segments $\Phi_1([0,p])$ and $\Phi_1([\tau /2,{(}1+\tau{)} /2])$ it foots on, hence the claim is true near the end points of the curve. Observe that the tangent vector to this simple curve is the unit vector  
$\left(\vartheta(z-p)/\vartheta(z+p)\right)^{a} $. 
It will be sufficient to show that this vector turns counter clockwise, i.e. that
\[
h(z,p)=\frac{\vartheta'(z-p)}{\vartheta(z-p)}-\frac{\vartheta'(z+p)}{\vartheta(z+p)}
\]
has constant sign along the curve. 
To do so, we write the $\vartheta$-function in terms of classical Weierstrass functions and use their well-known mapping properties: First, we have
\[
\frac{\vartheta'(z)}{\vartheta(z)}=\zeta(z) - \eta_1 z
\]
where $\zeta(z)$ is the Weierstrass $\zeta$-function and $\eta_1 = \int_{\alpha_1} \wp(z)\, dz$. Here $\wp(z)$ is the Weierstrass $\wp$-function, and $\alpha_1$ a cycle on $T$ homologous to $[0,1]$.

Secondly, the expression on the right hand side is known to map the rectangle $[0,1/2]\times[0,\tau /2]$ to the (suitably slit) right half plane (\cite{kob1}). In particular, the values of ${\vartheta'(z)}/{\vartheta(z)}$ will have positive real part. By symmetry and for any $p\in[0,1/2]$ and any $z\in [0,\tau/2]$, $h(z,p)$ will thus be a positive real number, as claimed.

The second case $p\to 0$ follows from the first by exchanging $p$ and $q$ and noting that the only thing that changes is the orientation of the integration path.
\end{proof}

\def\fw{2.2in}
\begin{figure}[H]
 \begin{center}
 \subfigure[(2,3,6)]{\label{2(2,3,6)}\includegraphics[width=2.0in]{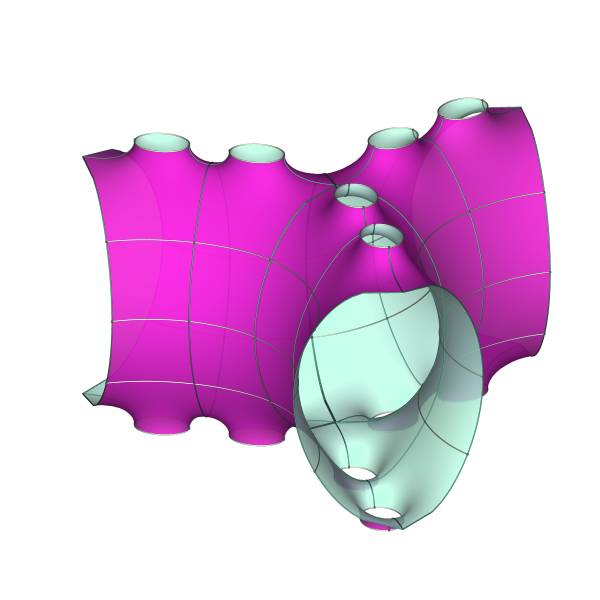}}
  \subfigure[(2,4,4)]{\label{2(2,4,4)}\includegraphics[width=\fw]{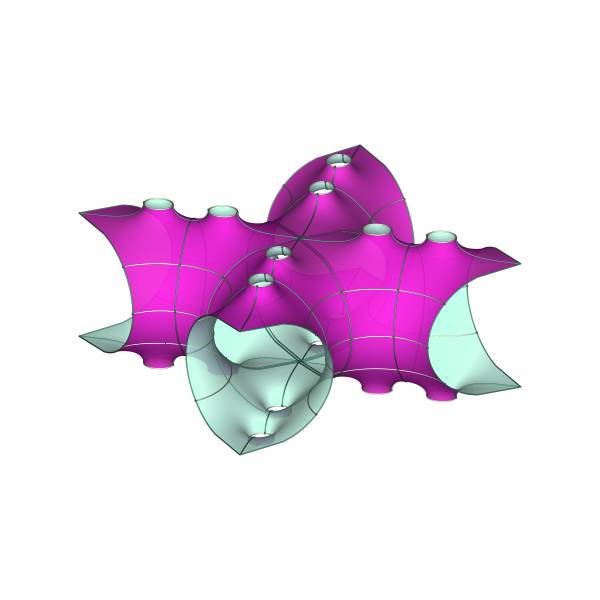}}\\
   \subfigure[(2,6,3)]{\label{2(2,6,3)}\includegraphics[width=\fw]{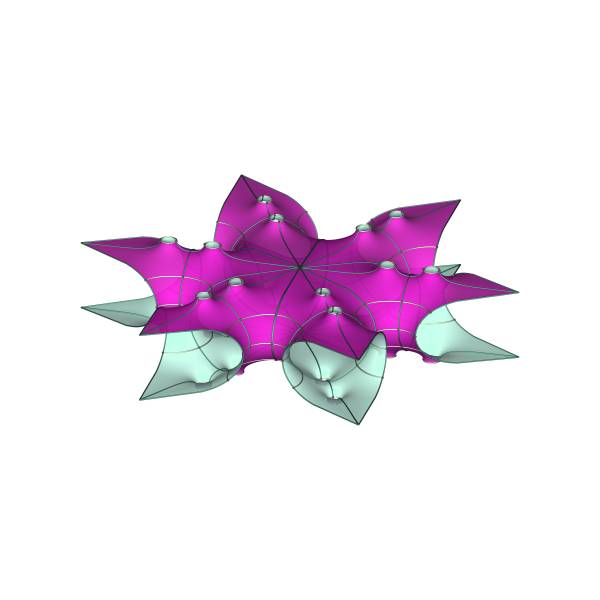}}
   \subfigure[(3,6,2)]{\label{2(3,6,2)}\includegraphics[width=\fw]{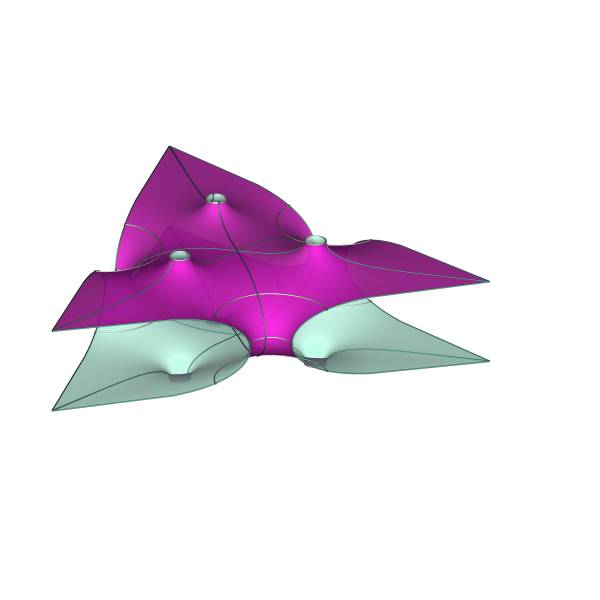}}
 \end{center}
 \caption{Symmetric case: equal exponents}
 \label{fig:sym1}
\end{figure}

The following table lists all possibilities of admissible angle combinations, except the cases of equal angles (i.e. $r=s=3$ or $r=s=4$). These cases reduce to the $(r,2,t)$ cases of   section \ref{sec:basic}, because the two hinges become symmetric, each being part of a $(r,2,t)$-triangle. 

%\begin{center}
%\begin{tabular}{l|ccccc}
%Name & $r$ & $s$ & $p$ vs $q$ & $\alpha_1(p)$ & $\alpha_2(p)=\pi/t$? \\ \hline
%Schoen H'-T-a & 2 & 3 & $6p+ 8q=1$ & $ -180$ & $30$ \\ 
%Schwarz P+h   & 2 & 4 & $4p+ 6q=1$ & $ -180$ & $45$ \\
%Lastcase      & 3 & 6 & $8p+10q=3$ & $-180$ & $90$ \\
%Schoen T-R'+h & 6 & 2 & $5p+ 3q=1$ & $  -180$ & $60$ 
%\end{tabular}
%\end{center}
%
\begin{center}
\begin{tabular}{cc|c}
$r$ & $s$ & Relation between $p$ and $q$ \\ \hline
2 & 3 & $6p+ 8\re(q)=1$ \\ 
2 & 4 & $4p+ 6\re(q)=1$ \\
 3 & 6 & $8p+10\re(q)=3$ \\
 6 & 2 & $5p+ 3\re(q)=1$ 
\end{tabular}
\end{center}

\subsection{Exponents have opposite signs}
\label{sec:opposite}

Without loss of generality, we can assume that $a<0$ and $b>0$.

In this case, the two hinges arrange themselves above each other and need to be matched up so that two edges are parallel and of the same length:

\def\fw{2in}
\begin{figure}[H]
 \begin{center}
   \subfigure{\includegraphics[width=\fw]{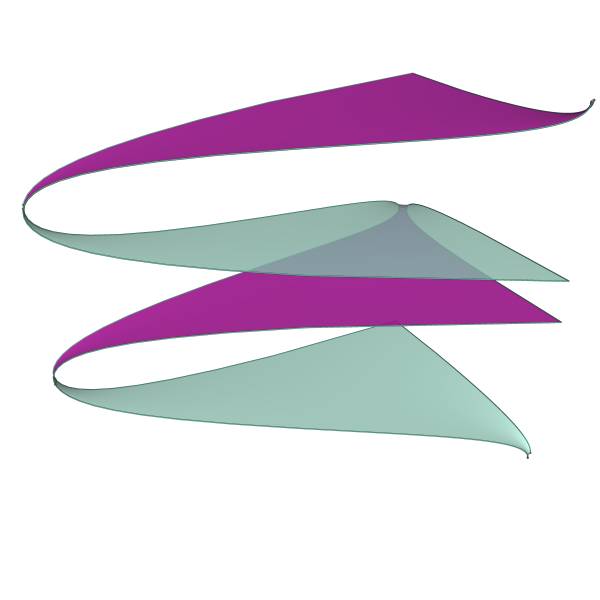}}
   \subfigure{\includegraphics[width=\fw]{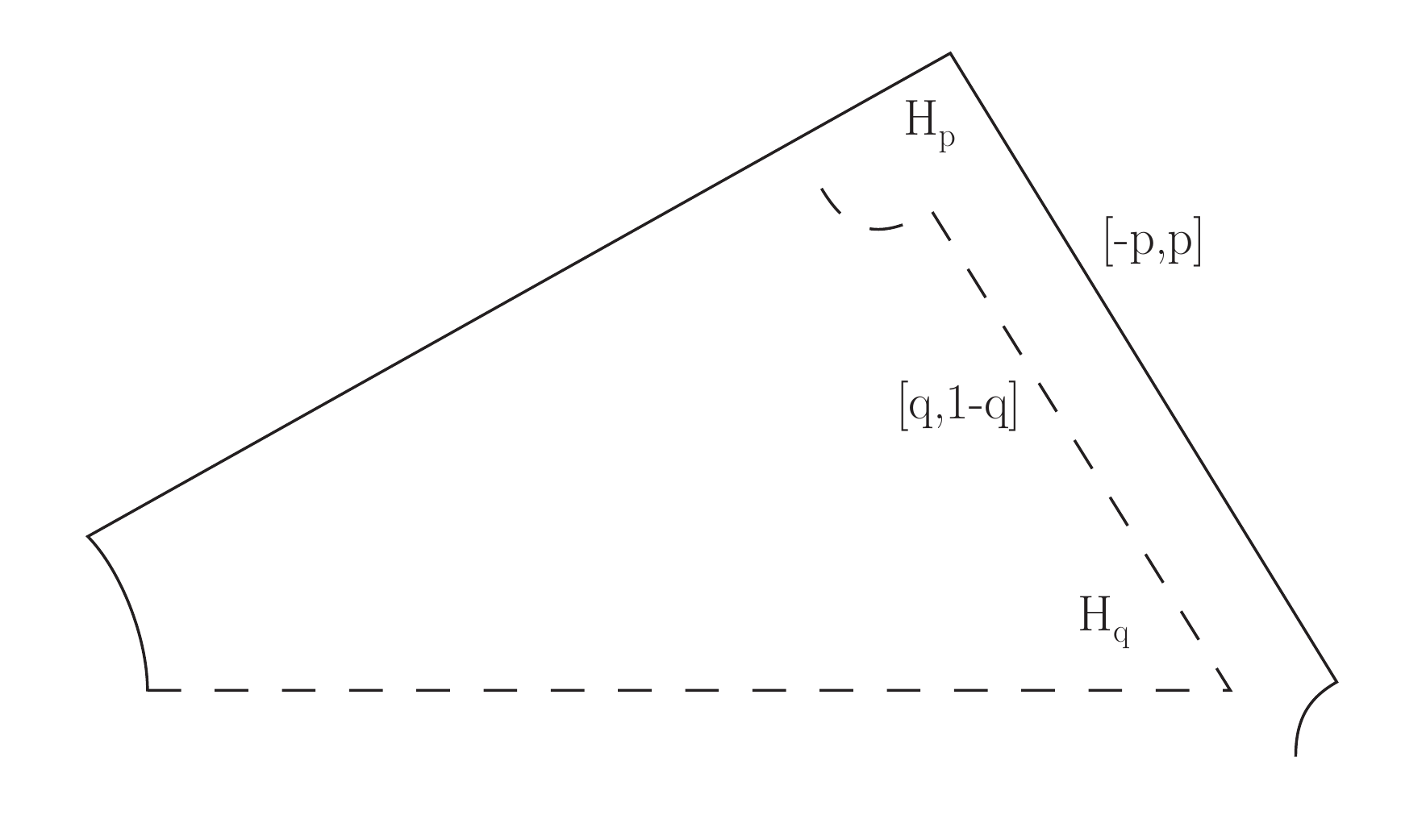}}
     \end{center}
 \caption{Minimal $(2,3,6)$-patch with schematic hinges of the symmetry planes }
 \label{fig:hinges2}
\end{figure}

\begin{figure}[h]
\centering
  \includegraphics[width=2in]{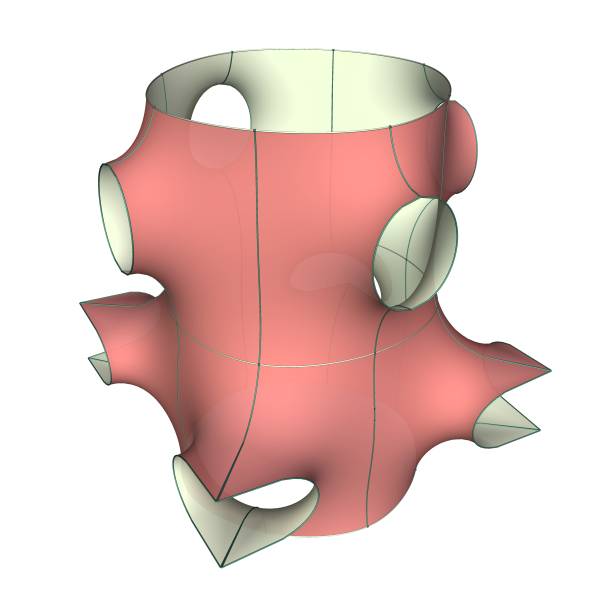} 
  \caption{Joining the basic cases $(2,4,4)$ and $(4,4,2)$ }
  \label{fig:basicjoin}
\end{figure}

There is a more intuitive way to understand the resulting surfaces geometrically, namely as vertically stacking two of the basic cases with the same underlying triangle group on top of each other. For instance, the $(2,4,4)$-surface (aka Schwarz P) and the $(4,4,2)$-surface (aka Schoen S'-S'') from the basic case family can be combined into the $(2,4,4)$-surface of the current family, see figure \ref{fig:basicjoin}. The different appearance compared to figure \ref{fig:sym2} is due to a different assembly of the surface from the fundamental piece $\Sigma$.

We normalize the divisors such that the edges $[-p,p]$ and $[q,{1}-q]$ are to be matched. 
Let $a = -(r-1)/r$ and $b = (s-1)/s$

This implies that we need the angle  $\alpha_0$ between $\Pi_{[-p,p]}$ and $\Pi_{[-q,q]}$ to be equal to $\pi/s$.

By proposition \ref{prop:bottom2}, the angle between the planes $\Pi_{[-p,p]}$ and $\Pi_{[-q,q]}$ is equal to $\pi \left( a(2p-1)+ b(2\re(q)-1)\right)$, and we obtain the constraint
\[
a(2p-1)+ b(2\re(q)-1)=\frac{1}{s}.
\]

Again, all combinatorially possible cases are listed in the table below. Two of the cases have been discovered earlier by A.~Schoen (I-WP) and H.~Karcher (T-WP). These are particularly simple in that $r=s$. This allows for more symmetric solutions to the period problem with $-p+\re(q)=1/(2(r-1))$. This results in horizontal straight lines on the surfaces, namely as images of the vertical lines through $p,q$ and $-p,-q$. Thus, the period problem is solved automatically. The other surfaces are probably new.

\begin{center}
\begin{tabular}{l|ccccc}
Name & $r$ & $s$ & Relation between $p$ and $q$ \\ \hline
            & 2 & 4 & $ 2p- 3\re(q)= -1$ &  \\
            & 2 & 6 & $ 6p-10\re(q)=-3$ &  \\ 
            & 2 & 3 & $ 6p- 8\re(q)=-3$ &  \\ 
            & 3 & 6 & $ 4p- 5\re(q)=-1$ &  \\
Schoen  I-WP & 4 & 4 & $ 6p= 3\re(q)= 1$ &  \\
Karcher T-WP & 3 & 3 & $24p= 8\re(q)= 3$ & 
\end{tabular}
\end{center}

\def\fw{2.2in}
\begin{figure}[H]
 \begin{center}
 \subfigure[(2,3,6)]{\label{3(2,3,6)}\includegraphics[width=\fw]{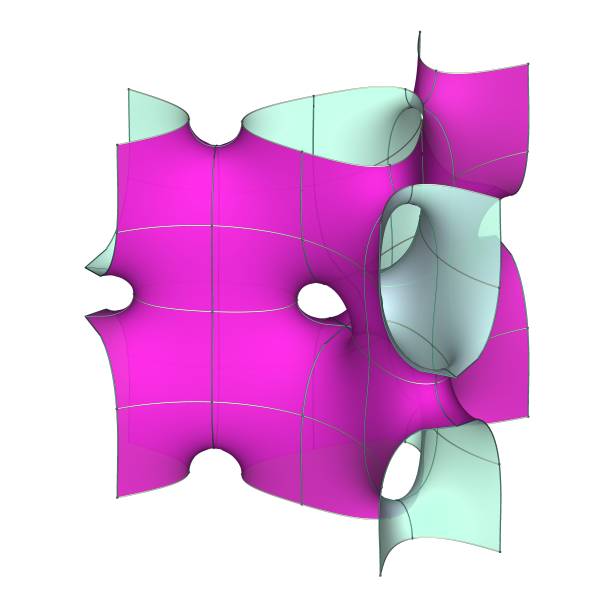}}
  \subfigure[(2,4,4)]{\label{3(2,4,4)}\includegraphics[width=\fw]{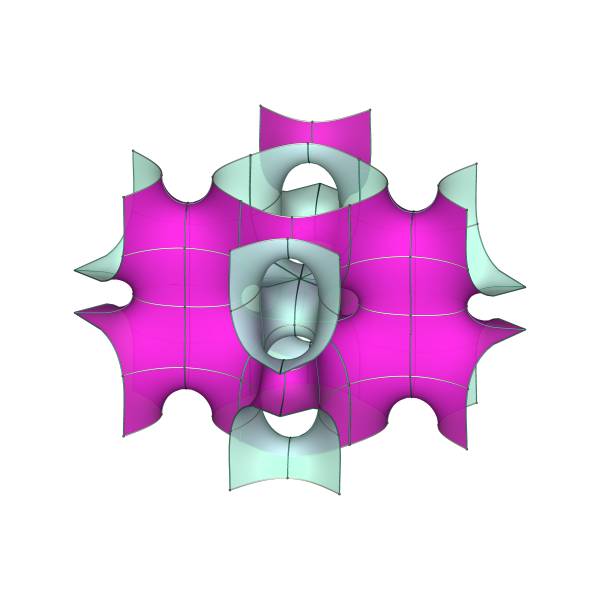}}\\
   \subfigure[(2,6,3)]{\label{3(2,6,3)}\includegraphics[width=\fw]{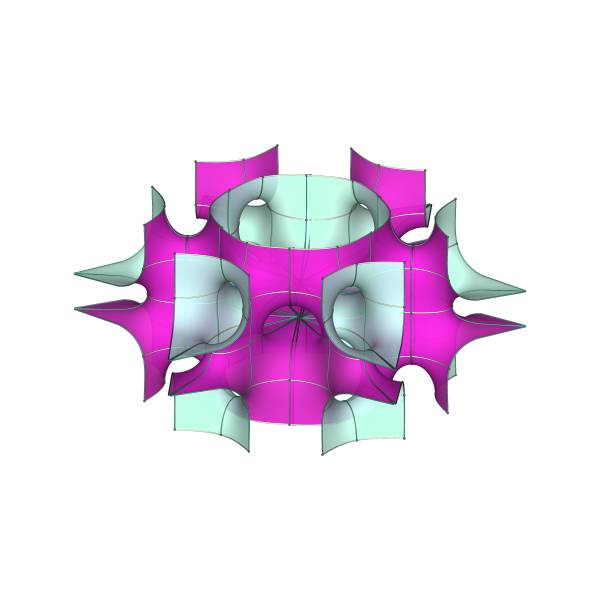}}
   \subfigure[(3,6,2)]{\label{3(3,6,2)}\includegraphics[width=\fw]{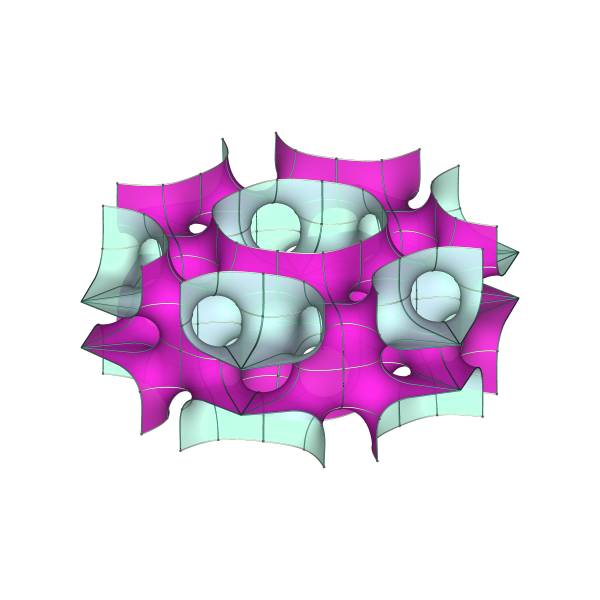}}\\
    \subfigure[(3,3,3)]{\label{3(3,3,3)}\includegraphics[width=\fw]{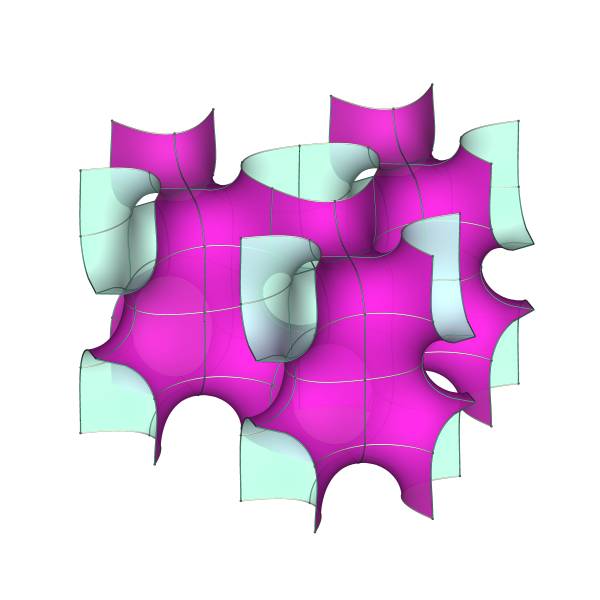}}
   \subfigure[(4,4,2)]{\label{3(4,4,2)}\includegraphics[width=\fw]{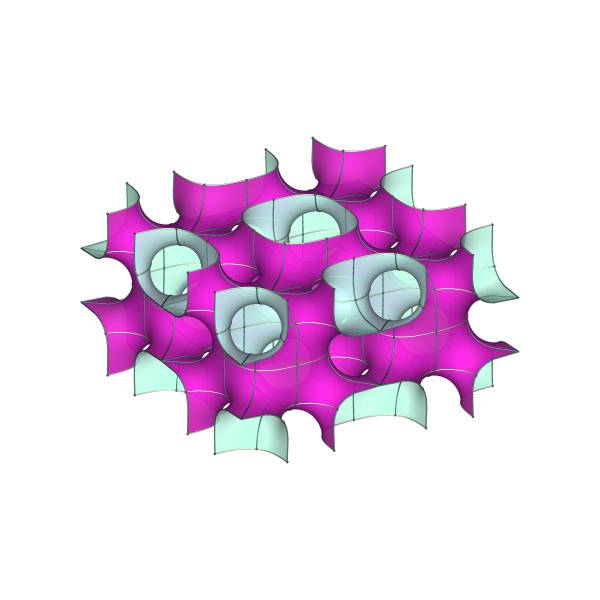}}

 \end{center}
 \caption{Symmetric case: opposite exponents}
 \label{fig:sym2}
\end{figure}

Note that any solution pair $(p,q)$ to the period problem is isometric to  $(1/2-p, 1/2-q)$ by shifting the divisor by $1/2$ and taking the reciprocal. The constraint equation will be slightly different then.

%\[
%  G = \frac{\vartheta (z+a)^{(m_1-1)/m_1}\vartheta (z-b-\tau)^{(m_2-1)/m_2}}
%           {\vartheta (z-a)^{(m_1-1)/m_1}\vartheta (z+b-\tau)^{(m_2-1)/m_2}}
%  \qquad \text{and} \qquad dh=dz
%\]

We will now show how to solve the period problem, using 
an extremal length argument.

\begin{theorem}
For any ${-1}<a<{0}$ and ${0}<b<{1}$, there is a  1-parameter family of values of any $\tau\in i\R^+$ such that the union of the hinges $H_p$ and $H_q$ forms a triangle. 
\end{theorem}
\begin{proof}

Recall that we labeled the points in $Z$ so that $p_{0}=-p$, $p_1=p$ and $q_{0}=-q$, $q_1={q}$.
The images of these points under a Schwarz-Christoffel map are denoted by $P_i$ and $Q_j$, respectively. 

Denote by $\Gamma_1$ the cycle in a periodic polygon  connecting $[P_0,P_1]$ with $[Q_{-1},Q_0]$, by $\Gamma_2$ the cycle connecting $[P_0,P_1]$ with $[Q_{1},Q_2]$. 

The period condition requires that all edges of $\Sigma$ corresponding to $[P_0,P_1]$, $[Q_{-1},Q_0]$ and $[Q_{1},Q_2]$ {lie in a same plane}. 
We will now reinterpret this condition in terms of the geometry of the periodic polygons $\Phi_1(Z)$ and $\Phi_2(Z)$. The periodic polygon $\Phi_1(Z)$ will have angles $\pi/r$ at $P_1$ and $\pi/s$ at $Q_0$. 
By rotating $G$, we can assume that the segments $[P_0,P_1]$ in both periodic polygons are on the real axes and point to the right.

By Proposition \ref{prop:symmetry}, the polygons  $\Phi_1(Z)$ and  $\Phi_2(Z)$ are symmetric with respect to a reflection about a vertical line, see figure \ref{fig:extlength}. Note, however, how the labeling of the corners changes.

To have $[P_0,P_1]$ to vertically line up with  $[Q_{-1},Q_0]$ in $\Sigma$, we need the imaginary parts of the periods of $\Gamma_1$ in $\Phi_1(Z)$ and $\Phi_2(Z)$ to be complex conjugate.

Equivalently, the period problem is solved if and only if the  horizontal segment from $Q_{-1}$ to $Q_0$ is at the middle of the height of the points $P_{-1}$ and $P_0$.

\begin{figure}[H] %  figure placement: here, top, bottom, or page
  \centering
  \includegraphics[width=4in]{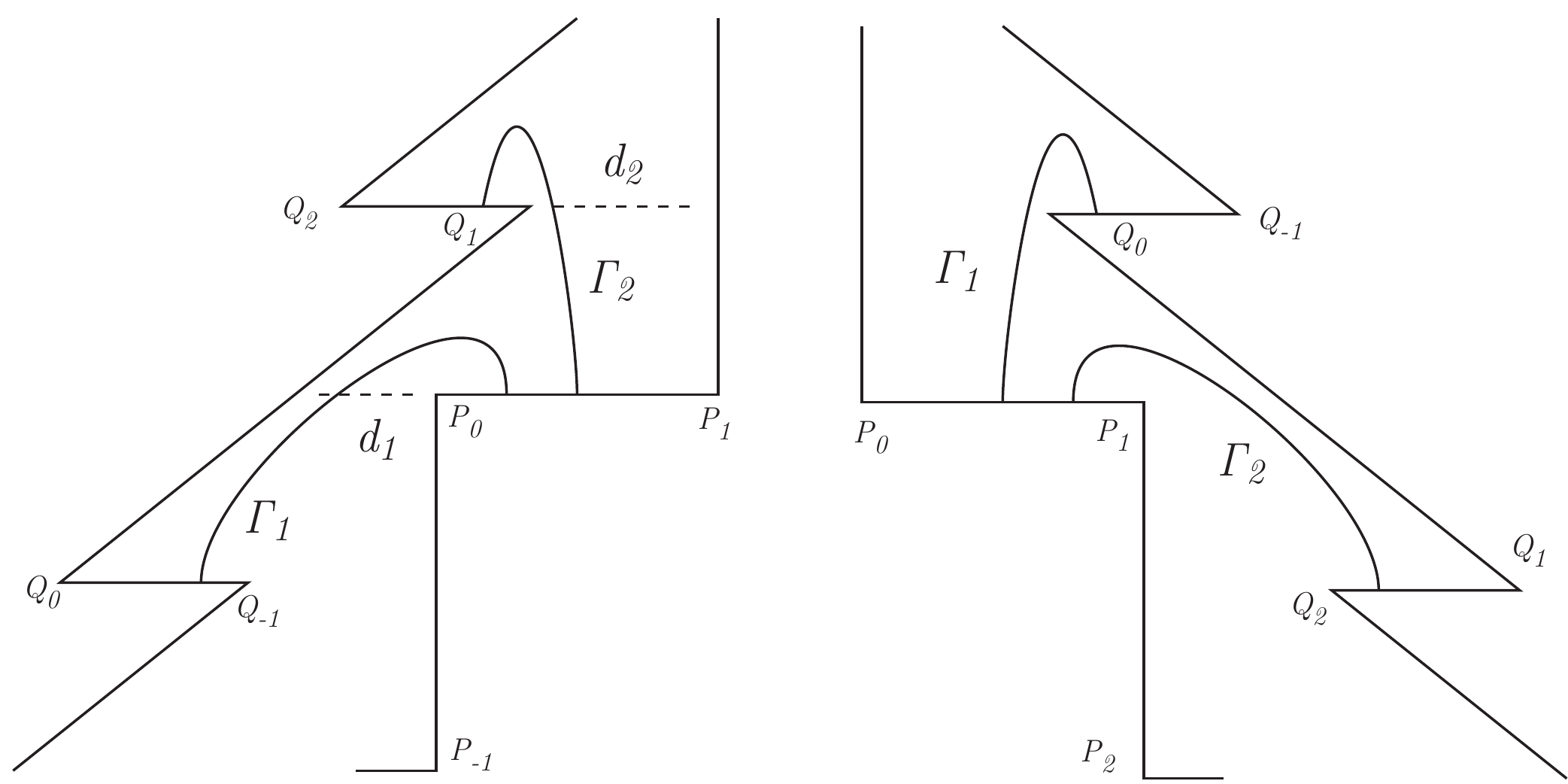} 
  \caption{Period condition and extremal length comparison}
  \label{fig:extlength}
\end{figure}

A periodic polygon with these properties can be uniquely constructed using as parameters a pair of numbers $d_1>0$ and $0<d_2<1$. Here $d_1$ measures the horizontal distance from $P_0$ to the opposite boundary arc, and $d_2$ the horizontal distance from $Q_1$ to the opposite boundary arc, see figure \ref{fig:extlength}.

However, such a periodic polygon does not necessarily correspond to a rectangular torus with points $p_0,p_1$ and $q_0,q_1$ placed symmetrically. By a translation, we can assume that $p_1=-p_0$ are symmetric, but not so for $q_0$ and $q_1$, see figure \ref{fig:exttorus}.

\begin{figure}[H]
\centering
  \includegraphics[width=3in]{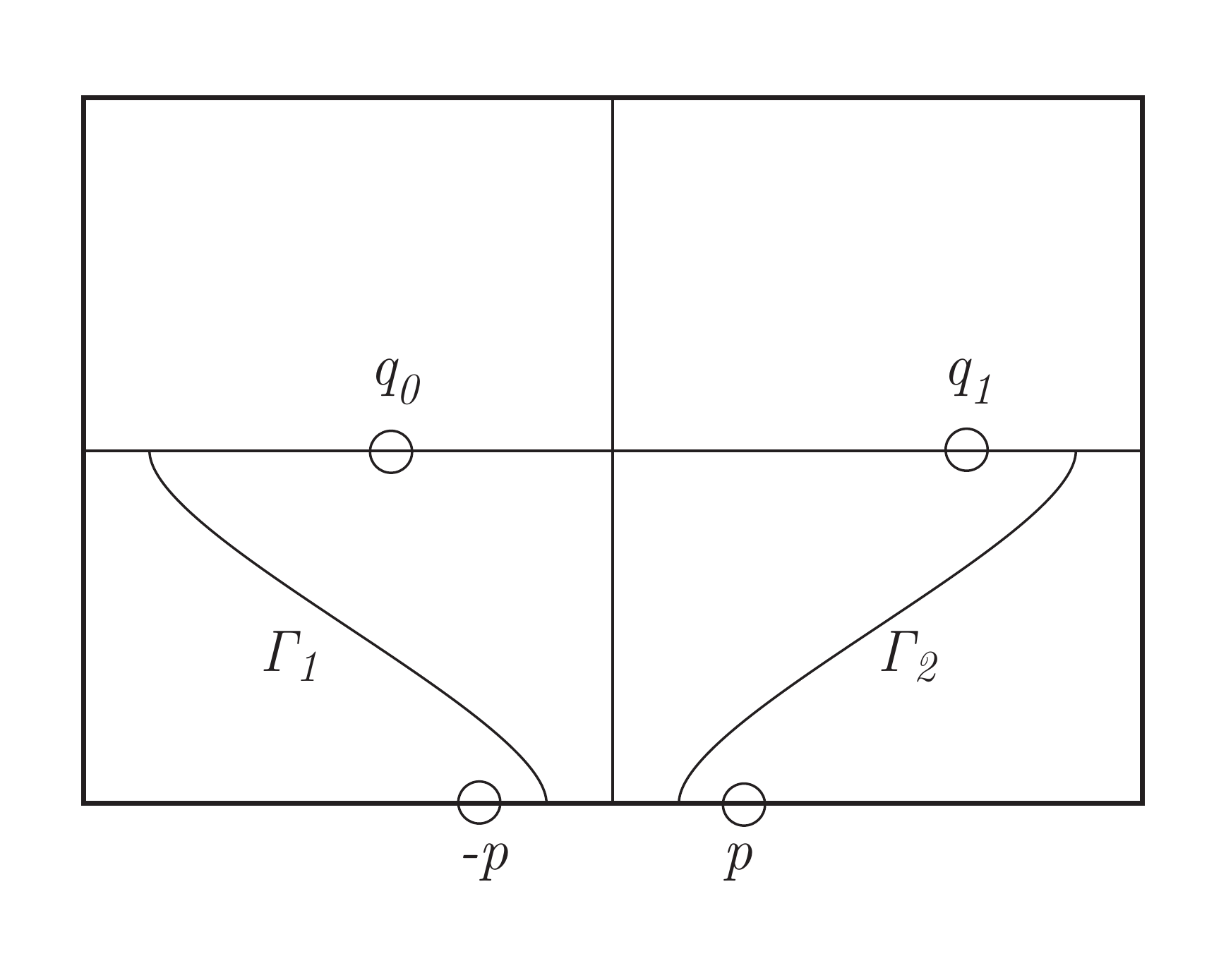} 
  \caption{Extremal length comparison}
  \label{fig:exttorus}
\end{figure}

Thus the period problem in this setting becomes the problem to find $d_1, d_2$ so that $Z$ enjoys the additional symmetry $q_1=-q_0$.

To measure the conformal symmetry of the periodic polygon, we use the extremal lengths $\ext(\Gamma_1)$ and $\ext(\Gamma_2)$ of the two cycles $\Gamma_1$ and $\Gamma_2$.  Then $q_1=-q_0$ if and only if $\ext(\Gamma_1)=\ext(\Gamma_2)$. This is because we can map $Z$ conformally to the upper half plane, where the cross ratio of the image points of $-p,p,q_{-1},q_0$ is determined by $\ext(\Gamma_1)$.

We will now construct a 1-parameter family of values for $(d_1,d_2)$ such that the corresponding periodic polygon is conformally symmetric as desired. To this end, we introduce a technical parameter $M$ for this family:

Let $f:(0,1)\to(0,\infty)$ be a bijection, and consider for fixed $M\in (0,\infty)$ the set $\Delta_M$ of all $(d_1,d_2)$ with $d_1+f(d_2)=M$. We  claim that for any $M>0$ there is at least one pair $(d_1,d_2)\in \Delta_M$  with $\ext(\Gamma_1)=\ext(\Gamma_2)$. 

If $d_1\to 0$, $d_2$ stays bounded away from $0$, and we obtain $\ext(\Gamma_1)\to \infty$ while $\ext(\Gamma_2)$ stays bounded. Similarly, when $d_2\to 0$, $\ext(\Gamma_2)\to \infty$ while $\ext(\Gamma_1)$ stays bounded. The claim now follows from the intermediate value theorem. Hence we obtain a 1-parameter family of solution of the period problem, with parameter $M$.
\end{proof}

Numerical evidence suggests that one can take $\tau\in i \R^+$ as the family parameter as well.

\section{Higher Genus}

The method above can be used to create many further examples of triply periodic minimal surfaces, by adding corners to the periodic polygons. This, however, also increases the dimension of the period problem. While there are methods available to solve such problems, it doesn't appear to be worth the effort at this point.

We briefly discuss two more cases which give examples of higher genus surfaces, which are related to surfaces that have been discussed in the literature.

\subsection{The Neovius Family}

For this family, we assume the reflectional family about the $x_1x_2$-plane as usual. The lower edge of $T$ is to have two corners at $\pm p$, while we assume four corners in the upper edge at $\pm q_1 +\tau/2$ and $\pm q_2 +\tau/2$.

We label the exponents at $p,q_1,q_2$ as $a,b_1,b_2$.

Our Neovius family has the exponents defined by

\begin{align*} 
  a  = {}& \frac1r-1<0,\\
  b_1= {}& 1-\frac1s>0,\\
  b_2= {}& 1-\frac1t>0.
\end{align*}

In the case $(r,s,t)=(2,4,4)$, we obtain a surface discovered by Schwarz' student Neovius (in the case of full cubic symmetry, \cite{ne1}). This case as well as the $(3,3,3)$ case 
(\cite{huff1}) allows for an additional symmetry with $p=1/4$ and $q_1+q_2=1/2$, which  renders the period problem 1-dimensional.

All other cases require to solve a 2-dimensional period problem,
and lead to 1-parameter families, whose existence we have established numerically.

\def\fw{2.2in}
\begin{figure}[H]
 \begin{center}
 \subfigure[(2,3,6)]{\label{n(2,3,6)}\includegraphics[width=\fw]{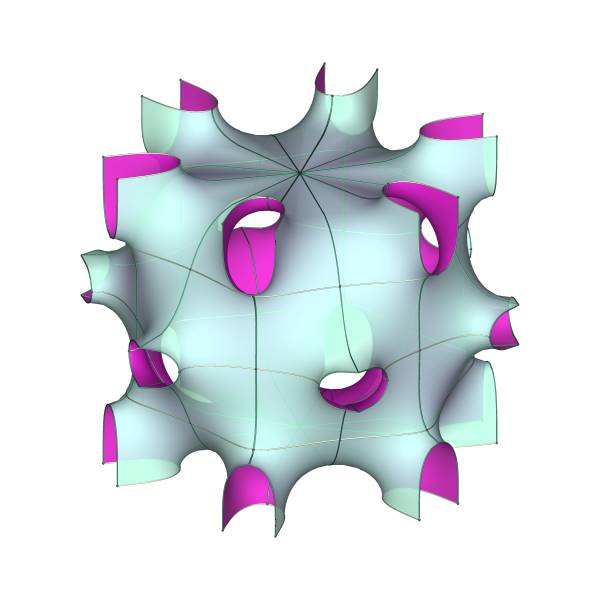}}
  \subfigure[(2,4,4)]{\label{n(2,4,4)}\includegraphics[width=\fw]{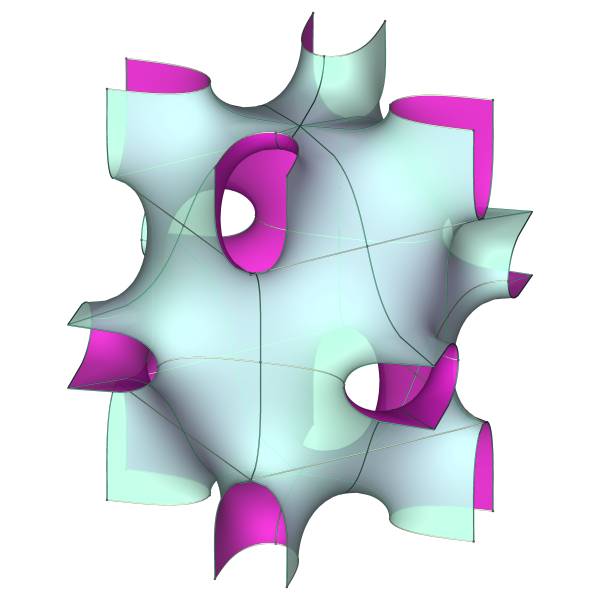}}\\
   \subfigure[(3,2,6)]{\label{n(3,2,6)}\includegraphics[width=\fw]{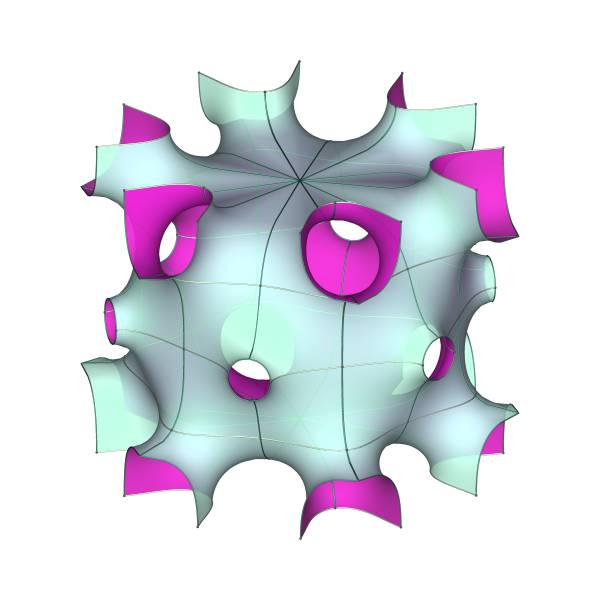}}
   \subfigure[(3,3,3)]{\label{n(3,3,3)}\includegraphics[width=\fw]{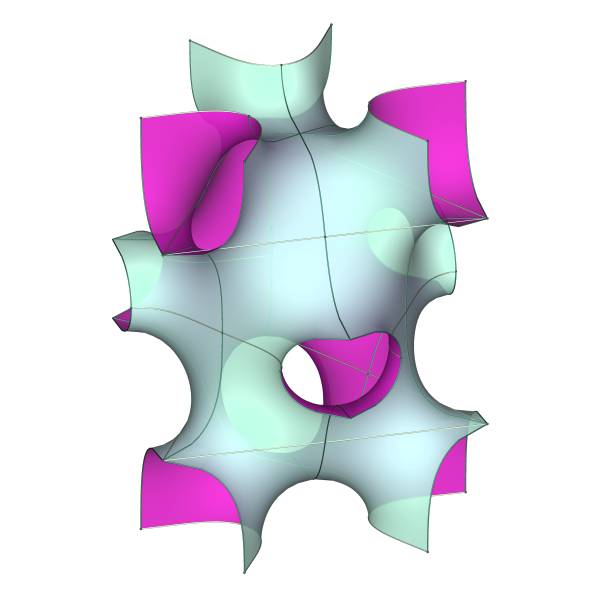}}\\
    \subfigure[(4,2,4)]{\label{n(4,2,4)}\includegraphics[width=\fw]{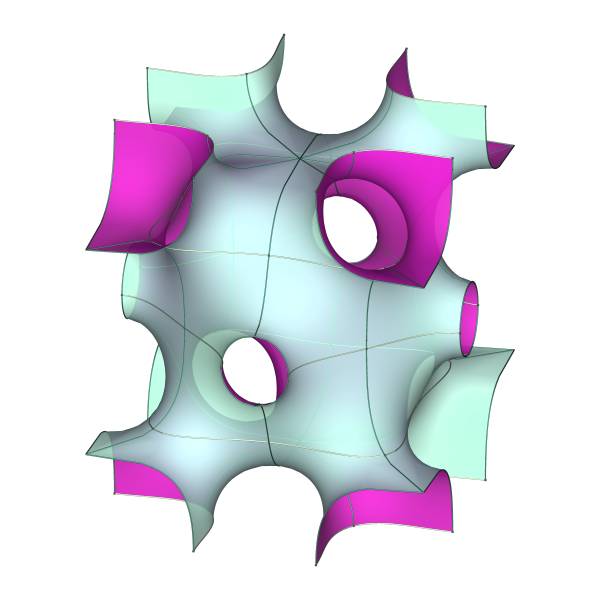}}
   \subfigure[(6,2,3)]{\label{n(6,2,3)}\includegraphics[width=\fw]{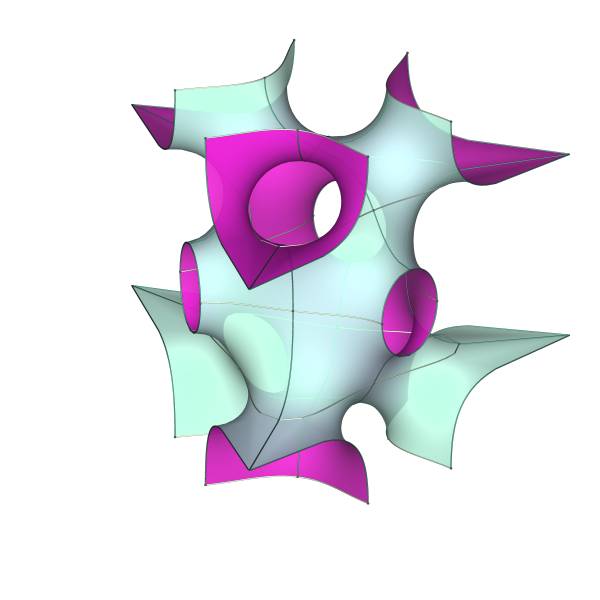}}

 \end{center}
 \caption{Neovius Family}
 \label{fig:neovius}
\end{figure}

\subsection{The Multiple Spout Families}

One of the motivations behind this paper was to understand the possibility of creating minimal vertical cylinders with spout-like openings pointing in several directions, where the ``spouts'' are bounded by planar symmetry curves meeting at the same angle for each direction. By choosing the angles suitably and closing the periods so that the tips of the spouts line up, replicating such a surface by reflection gives embedded triply periodic minimal surfaces.

\def\fw{1.5in}
\begin{figure}[H]
 \begin{center}
 \subfigure[]{\includegraphics[width=\fw]{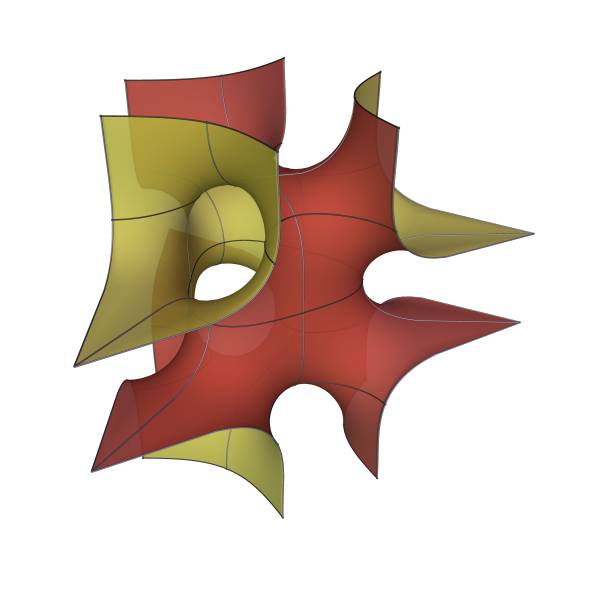}}
  \subfigure[]{\includegraphics[width=\fw]{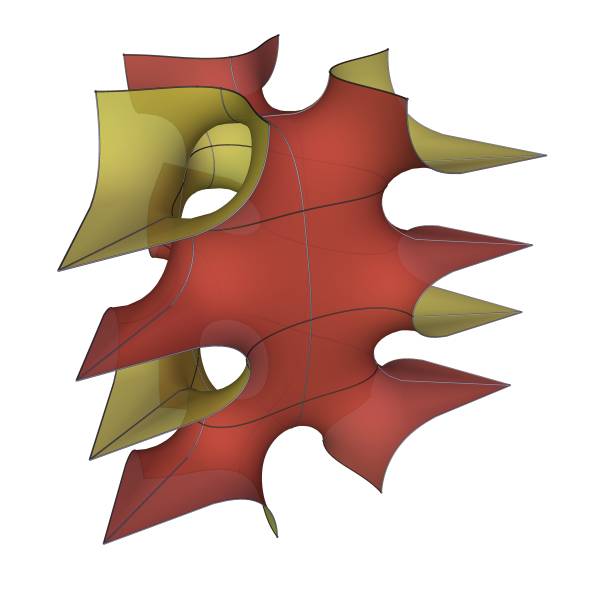}}
   \subfigure[]{\includegraphics[width=\fw]{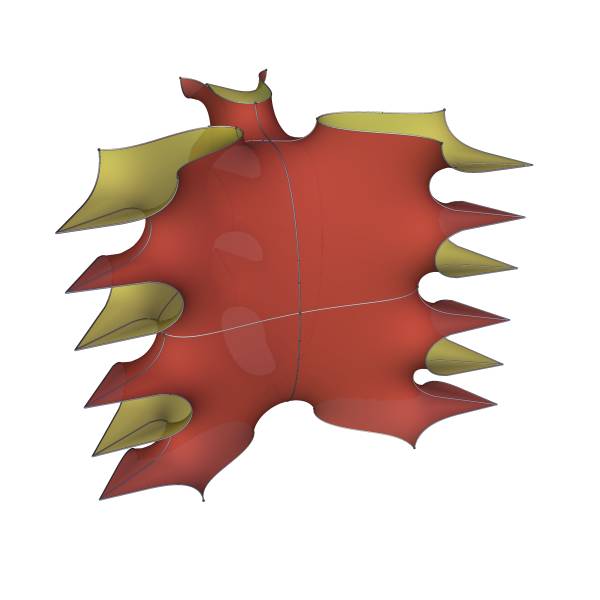}}
    \end{center}
 \caption{Minimal cylinders with spout-like openings}
 \label{fig:spouts}
\end{figure}

Again we assume the reflectional family about the $x_1x_2$-plane as usual. The lower edge of $T$ is to have to corners at $\pm p$, while we assume $2n$ corners in the upper edge at $\pm q_i +\tau/2$, $i=1,\ldots,n$.

We label the exponents at $p$  {(}resp. $q_i${)} as $a$ {(}resp. $b_i${)} and set $a_1=(1-r)/r<0$ and $b_i= (-1)^i(s-1)/s$.

This results in surfaces  with $r$ directions into which the spouts point, and spout angles of $2\pi/s$.

Obtaining embedded surfaces requires to solve an $n$-dimensional period problem. In Figure \ref{fig:spouts}, we show solutions for $(r,s,t)=(3,6,2)$ for $n=1,2,3$.

The surfaces obtained in these families are closely related to the minimal surfaces obtained by Traizet in section 5.2 of \cite{tr3a}. 
There he describes surfaces that can be obtained by gluing singly periodic Scherk surfaces together, which have been placed at the vertices of a periodic tiling of the plane. However, his theorem does not apply to all cases, as the underlying tilings are not rigid (in Traizet's sense). 

Below are pictures for all combinatorially possible distinct cases for $n=2$.

\def\fw{1.6in}
\begin{figure}[H]
 \begin{center}
 \subfigure[(2,3,6)]{\label{4(2,3,6)}\includegraphics[width=\fw]{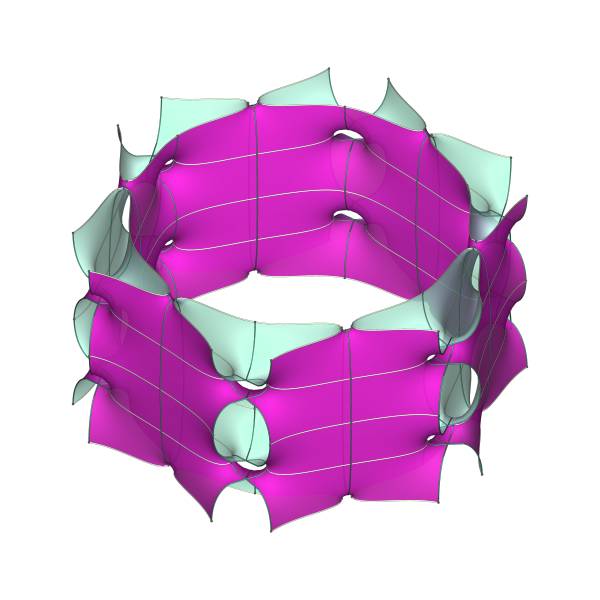}}
  \subfigure[(2,4,4)]{\label{4(2,4,4)}\includegraphics[width=\fw]{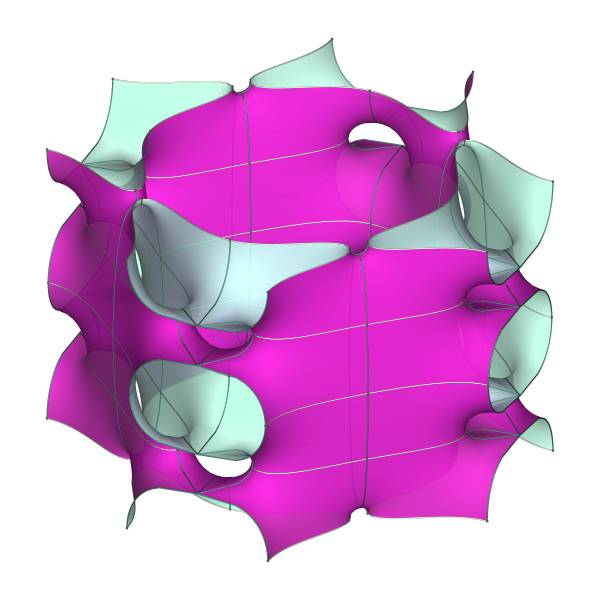}}
   \subfigure[(2,6,3)]{\label{4(2,6,3)}\includegraphics[width=\fw]{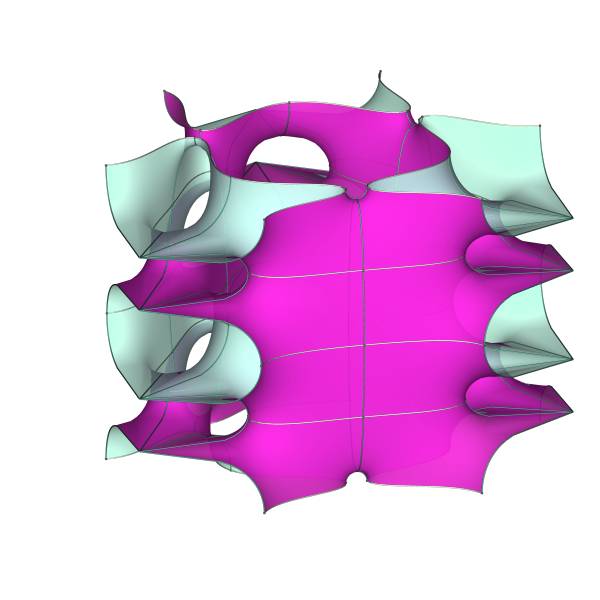}}\\
   \subfigure[(3,2,6)]{\label{4(3,2,6)}\includegraphics[width=\fw]{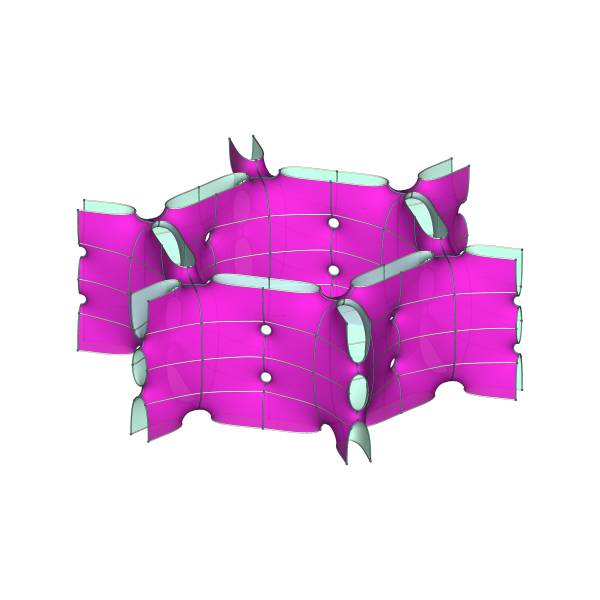}}
    \subfigure[(3,3,3)]{\label{4(3,3,3)}\includegraphics[width=\fw]{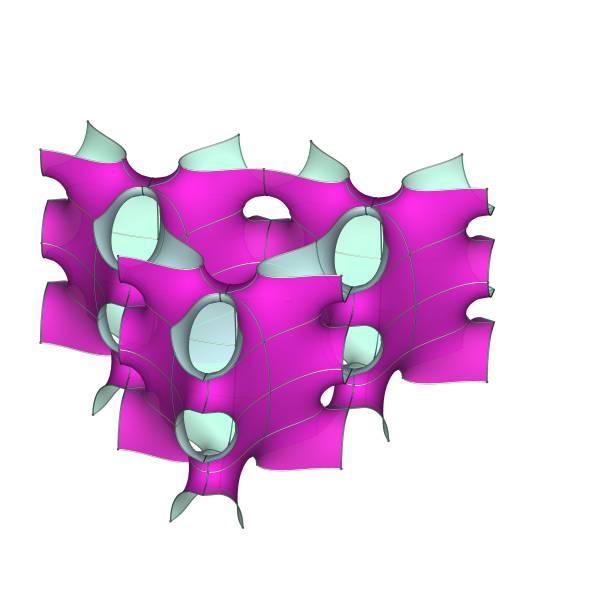}}
    \subfigure[(3,6,2)]{\label{4(3,6,2)}\includegraphics[width=\fw]{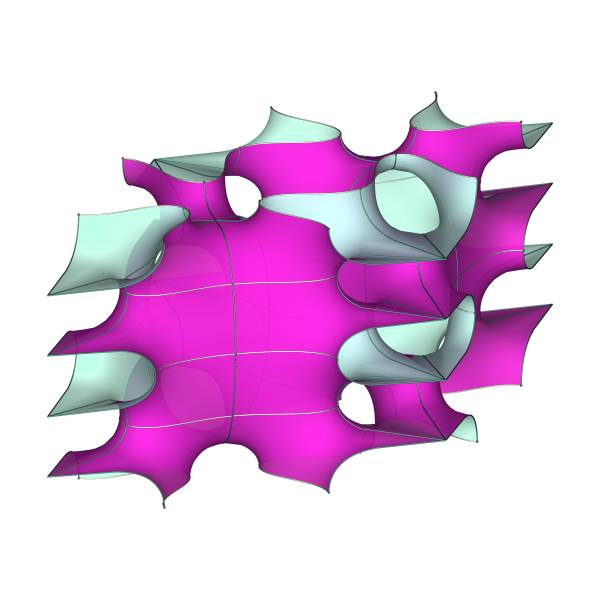}}\\
   \subfigure[(4,2,4)]{\label{4(4,2,4)}\includegraphics[width=\fw]{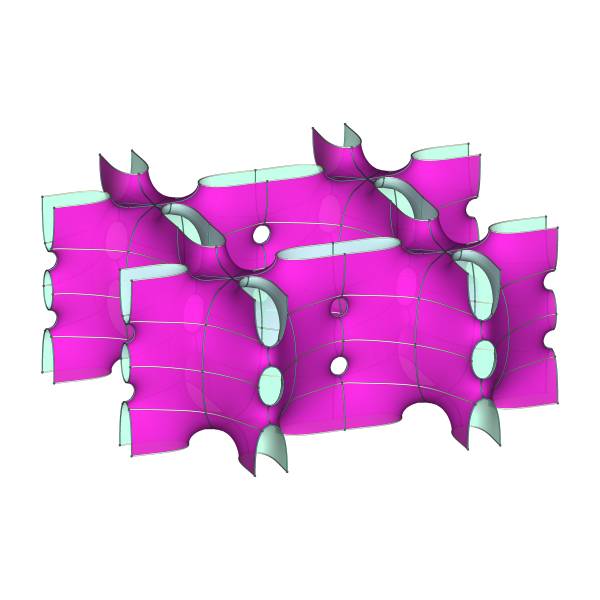}}
  \subfigure[(4,4,2)]{\label{4(4,4,2)}\includegraphics[width=\fw]{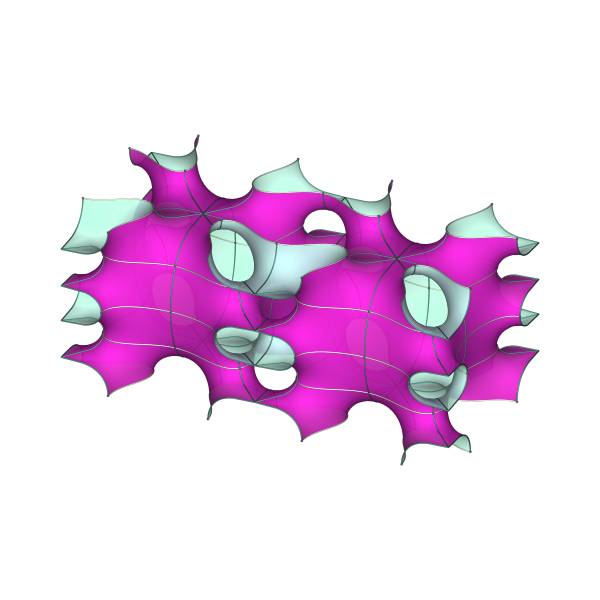}}
  \subfigure[(6,2,3)]{\label{4(6,2,3)}\includegraphics[width=\fw]{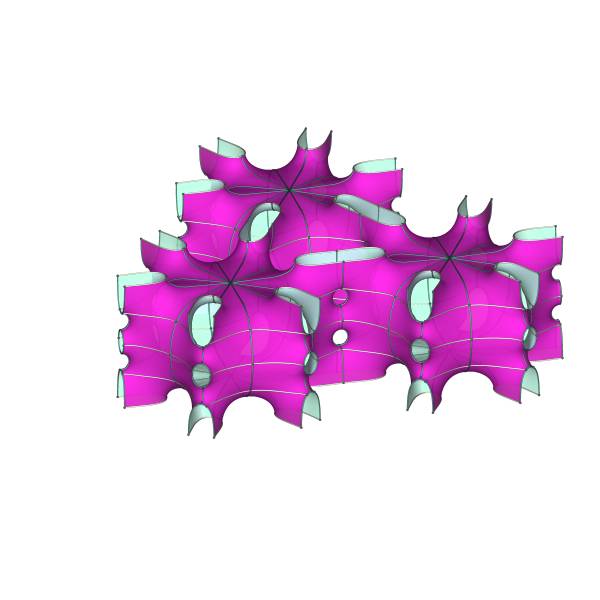}}\\
\subfigure[(6,3,2)]{\label{4(6,3,2)}\includegraphics[width=\fw]{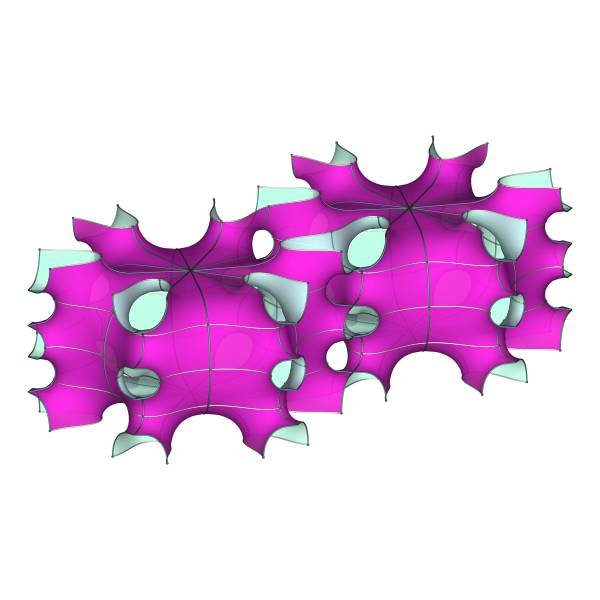}}

 \end{center}
 \caption{Double Spout Family}
 \label{fig:double}
\end{figure}

\section{Embeddedness}

In this section, we show {the following:}

\begin{theorem}
All of the triply periodic minimal surfaces discussed above are {indeed} embedded.
\end{theorem}
\begin{proof} The approach is quite standard, and goes as follows:

Instead of looking at $\Sigma$, we will look at the conjugate surface $\Sigma^*$. The piece corresponding to the domain $[0,1/2]\times[0,\tau/2]$ will be denoted by $\Sigma^*_0$. By our symmetry assumption, it has polygonal boundary.

According to Rad\'o and Nitsche (\cite{rad1,nit1}), a given Jordan curve in $\R^3$ which projects monotonely onto a convex curve in a plane,  has a unique Plateau solution which is a graph over the convex domain bounded by the planar curve.

Then, according to Krust (\cite{ka6}), the conjugate of $\Sigma^*_0$ will also be a graph, and in particular be embedded
within the fundamental prism determined by the reflection triangle and the horizontal symmetry planes.

Finally, reflecting at the symmetry planes will only generate disjoint copies, leaving the entire surface embedded.

Thus, all we have to show is that the polygonal boundary of $\Sigma^*_0$ is a graph over a convex domain, possibly except
for finitely many vertical segments. 

The polygonal boundary contour of $\Sigma^*_0$ consists of horizontal segments, each perpendicular to the symmetry plane in which the corresponding planar symmetry curve of $\Sigma_0$ lies. In addition, the image segments of the boundary intervals $[0,\tau/2]$ and $[1/2,(1+\tau)/2]$ are vertical segments. The angles between consecutive horizontal segments are the angles of the reflection group triangles.

We will now discuss the embeddedness of the surfaces
from the family in section \ref{sec:opposite}. From the angle condition we deduce that the two image segments of $[0,p]$ and $[q,(1+\tau)/2]$ are parallel. Thus we can assume without loss that they are parallel to the ${x_1}$-axis. In addition, the image segments of  $[0,\tau/2]$ and $[1/2,(1+\tau)/2]$ are parallel to the ${x_3}$-axis, as noted before. As the angles at the corners corresponding to $p$ and $q$ are $\pi/r$ and $\pi/s$, respectively, the entire boundary contour lies above a rectangle in the ${x_2x_3}$-plane, and is a graph except for the two vertical segments. 

The argument for the other families is even simpler.

\begin{figure}[H] %  figure placement: here, top, bottom, or page
  \centering
  \includegraphics[width=4in]{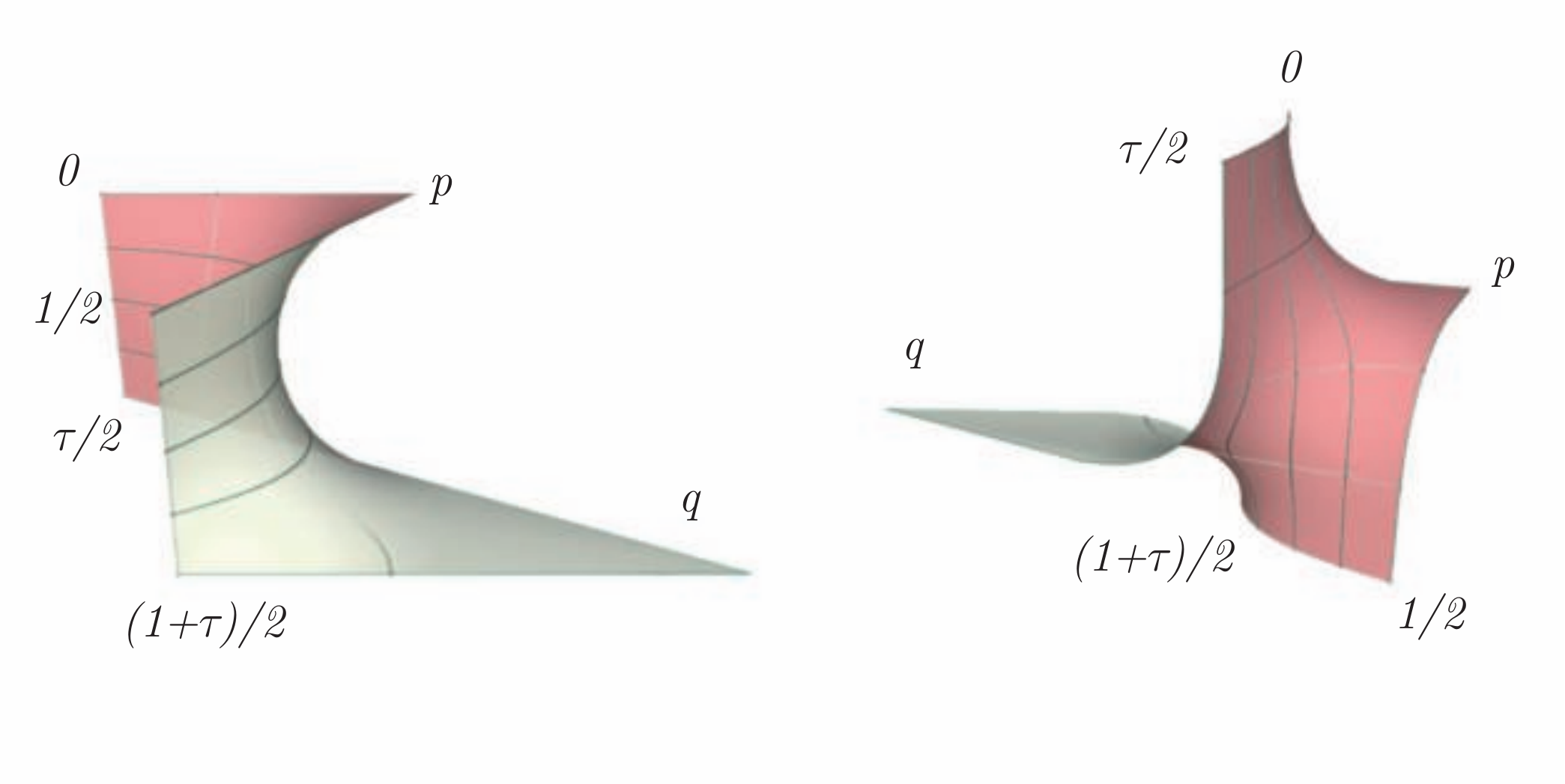} 
  \caption{Conjugate and original patch of the $(3,6,2)$-surface}
  \label{fig:embedded}
\end{figure}

\end{proof}

%\cleardoublepage
%\addcontentsline{toc}{section}{References}

\bibliographystyle{plain}
\bibliography{min}

\begin{thebibliography}{10}

\bibitem{fk1}
W.~Fischer and E.~Koch.
\newblock On 3-periodic minimal surfaces.
\newblock {\em Zeitschrift f\"{u}r Kristallographie}, 179:31--52, 1987.

\bibitem{fk2}
W.~Fischer and E.~Koch.
\newblock On 3-periodic minimal surfaces with non-cubic symmetry.
\newblock {\em Zeitschrift f\"{u}r Kristallographie}, 183:129--152, 1988.

\bibitem{huff1}
R.~Huff.
\newblock Flat structures and the triply periodic minimal surfaces ${C(H)}$ and
  ${tC(P)}$.
\newblock {\em Houston J. Math.}, 32:1011--1027, 2006.

\bibitem{ka6}
H.~Karcher.
\newblock Construction of minimal surfaces.
\newblock {\em Surveys in Geometry}, pages 1--96, 1989.
\newblock University of Tokyo, 1989, and Lecture Notes No. 12, SFB256, Bonn,
  1989.

\bibitem{ka5}
H.~Karcher.
\newblock The triply periodic minimal surfaces of alan schoen and their
  constant mean curvature companions.
\newblock {\em Manuscripta Math.}, 64:291--357, 1989.

\bibitem{kapo1}
H.~Karcher and K.~Polthier.
\newblock Construction of triply periodic minimal surfaces.
\newblock {\em Phil. Trans. R. Soc. Lond.}, pages 2077--2104, 1996.

\bibitem{kob1}
H.~Kober.
\newblock {\em Dictionary of Conformal Representation}.
\newblock Dover, 2nd edition, 1957.

\bibitem{me6}
W.~H. Meeks~III.
\newblock The theory of triply-periodic minimal surfaces.
\newblock {\em Indiana Univ. Math. J.}, 39(3):877--936, 1990.

\bibitem{mum1}
D.~Mumford.
\newblock {\em Lectures on Theta I}.
\newblock Birkh\"{a}user, Boston, 1983.

\bibitem{ne1}
E.~R. Neovius.
\newblock {\em Bestimmung zweier spezieller periodischer Minimalfl\"{a}chen}.
\newblock Akad. Abhandlungen, Helsingfors, 1883.

\bibitem{nit1}
J.~{C.}~{C.} Nitsche.
\newblock {\em Lectures on Minimal Surfaces}, volume~1.
\newblock Cambridge University Press, 1989.

\bibitem{rad1}
T.~Rado.
\newblock {\em On the problem of Plateau}.
\newblock Springer Verlag, Berlin, 1933.

\bibitem{Ra1}
V.~Ramos~Batista.
\newblock A family of triply periodic costa surfaces.
\newblock {\em Pacific J. Math.}, 212:347--370, 2003.

\bibitem{sch1}
A.~Schoen.
\newblock Infinite periodic minimal surfaces without self-intersections.
\newblock Technical Note D-5541, NASA, Cambridge, Mass., May 1970.

\bibitem{schw1}
H.~A. Schwarz.
\newblock {\em Gesammelte Mathematische Abhandlungen}, volume~1.
\newblock Springer, Berlin, 1890.

\bibitem{tr3a}
M.~Traizet.
\newblock Construction of triply periodic minimal surfaces.
\newblock preprint 172 Tours, 1998.

\bibitem{tr7}
M.~Traizet.
\newblock On the genus of triply periodic minimal surfaces.
\newblock {\em J. Differential Geometry}, 2007.
\newblock to appear.

\end{thebibliography}
%\bibliography{liter}

\bigskip
%%%%%%%%%%%% Authors addresses %%%%%%%%%%%%%
{\small
\noindent
Shoichi Fujimori \\
Department of Mathematics \\
Fukuoka University of Education \\
Munakata, Fukuoka 811-4192 \\
Japan \\
{\itshape E-mail address}\/: fujimori@fukuoka-edu.ac.jp
\par\vskip4ex
\noindent
Matthias Weber \\
Department of Mathematics \\
Indiana University \\
Bloomington, IN 47405 \\
USA \\
{\itshape E-mail address}\/: matweber@indiana.edu
}
\end{document}